\documentclass[11pt]{article}

\title{ On $\psi$-umbral extensions  of Stirling numbers and
Dobinski-like formulas}

\author{A. Krzysztof Kwa\'sniewski\\  
\\ High School of Mathematics and Applied Informatics\\
PL - 15-021 Bialystok , ul.Kamienna 17,  Poland
\\e-mail: kwandr@wp.pl}

\usepackage{amsmath,amsthm}

\chardef\bslash=`\\ 
\begin{document}
\maketitle\
\begin{abstract}

\noindent A so called $\psi$-umbral extensions of the  Stirling
numbers of the second kind are considered and the resulting
Dobinski-like various formulas - including \textit{new} ones - are
presented. These extensions naturally encompass the two well known
$q$-extensions. The further consecutive \textit{$\psi$- umbral
extensions} of Carlitz-Gould-Milne $q$-Stirling numbers are
therefore realized here in a two-fold way. The fact that the
umbral $q$-extended Dobinski formula may also be interpreted as
the average of powers of random variable $X_q$ with the
$q$-Poisson distribution singles out the $q$-extensions which
appear to be a kind of "singular point" in the domain of
$\psi$-umbral extensions as expressed by Observations $2.1$ and
$2.2$. Other relevant possibilities are tackled with the paper`s
closing down questions and suggestions with respect to other
already existing extensions while a brief limited survey of these
other type extensions is being delivered. There the Newton
interpolation formula and divided differences appear helpful and
inevitable along with umbra symbolic language in describing
properties of general exponential polynomials of Touchard and
their possible generalizations. Exponential structures or
algebraically  equivalent prefabs with their exponential formula
appear to be also naturally relevant.

\end{abstract}

\noindent MCS numbers:  05A40, 11B73, 81S99

\noindent formulas, Graves-Heisenberg-Weil algebra.

\vspace{0.5cm}

\subsection*{1. In the $q$-extensions realm}

\vspace{0.2cm}

At first let us make a remark on notation (see also Appendix).
$\psi$ is a number or functions` sequence - sequence of functions
of a parameter $q$. $\psi$ denotes an extension of
$\langle\frac{1}{n!}\rangle_{n\geq 0}$ sequence to quite arbitrary
one (the so called - "admissible" [1, 2]). The specific choices
are for example : Fibonomialy-extended sequence
$\langle\frac{1}{F_n!}\rangle_{n\geq 0}$ ($\langle F_n \rangle$ -
Fibonacci sequence ) or just "the usual" $\psi$-sequence
$\langle\frac{1}{n!}\rangle_{n\geq 0}$ or Gauss $q$-extended
$\langle\frac{1}{n_q!}\rangle_{n\geq 0}$ admissible sequence of
extended umbral operator calculus, where $ n_q=\frac{1-q^n}{1-q}$
and $n_q!=n_q(n-1)_q!, 0_q!=1$ - see more below. With such type
extension we may "$\psi$-mnemonic" repeat with exactly the same
simplicity this what was done by Rota forty one years ago. Namely
forty one years ago Gian-Carlo Rota \cite{3} proved that the
exponential generating function for Bell numbers $B_n$ is of the
form

\begin{equation}\label{eq1}
     \sum_{n=0}^\infty \frac {x^n}{n!}B_n= \exp(e^x-1)
\end{equation}
using the linear functional  \textit{L } such that
\begin{equation}\label{eq2}
     L(X^{\underline{n}})=1, \qquad  n\geq 0.
\end{equation}
Bell numbers (see: formula (4)in \cite{3}) are then defined by
\begin{equation}\label{eq3}
          L(X^n)=B_n,\qquad   n \geq 0.
\end{equation}
The above formula is exactly the Dobinski formula \cite{4} if $L$
is interpreted as the average value functional for the random
variable $X$ with the Poisson distribution where $L(X) = 1$. As a
matter of fact it is Blissard calculus inspired umbral formula
\cite{3} (see [5] for umbral nature of Poisson random variables
and the introduction in [5] for historical remarks on Blissard`s
calculus roots). On this occasion let us recall that the Stirling
numbers of the second kind are relatives of the Poisson
distribution in the known way. Namely if $X$ is a random variable
with a Poisson distribution with expected value $\lambda$, then
its $n-th$ moment is the exponential polynomial  $\varphi_n$ value
at $\lambda$  i.e.
$$E(X^n)=\varphi_n(\lambda)=\sum_{k=0}^{n}\Big\{ {n \atop k}\Big\}\lambda^k.$$
Hence in particular, the $n-th$ moment of the Poisson distribution
with expected value $1$ is precisely the number of partitions of a
set of size  $n$ i.e. it is the $n-th$ Bell number (this fact is
Dobinski's formula as stated by the formula (3)). The formula (3)
is tempting to be $\psi$-extended somehow as the $\psi$-Poisson
process distribution is known [2, 1]. Before doing this let us
remind that recently an interest to extensions of Stirling numbers
and consequently to Bell numbers was revived among "$q$-coherent
states physicists" \cite{6, 7, 8} with several important
generalizations already at hand such as in [9, 10, 11]. The merit
of such applications  is in that the expectation value with
respect to coherent state $|\gamma>$ with $|\gamma| = 1$ of the
$n$-th power of the number of quanta operator \cite{6, 7, 8} is
"just" the $n$-th Bell number $B_n$ and the explicit formula for
this expectation number of quanta is "just" Dobinski formula
\cite{6}. The same holds for $q$-coherent states case \cite{6}
i.e. the expectation value with respect to $q$-coherent state
$|\gamma> $ with $|\gamma| = 1$ of the $n$-th power of the number
operator is the $n$-th $q$-Bell number \cite{8, 6} defined as the
sum  of $q$-Stirling numbers  $ \Big\{{n \atop k}\Big\}_q $
introduced  by Carlitz and Gould and recently exploited among
others in \cite{12, 6, 7, 8}. Note there then that for the
\textbf{two standard }\cite{12} \textbf{$q$-extensions} of the
Stirling numbers of the second kind  we have as the \textbf{first}
ones the  $q$-Stirling numbers:

\begin{equation}\label{eq4}
x_q^n=\sum_{k=0}^{n}\Big\{{n \atop k}\Big\}_q  x_q^{\underline k},
\end{equation}
where $ x_q=\frac{1-q^x}{1-q}$ and $x_q^{\underline
k}=x_q(x-1)_q...(x-k+1)_q $  and then the\textbf{ second} ones
called $q^\sim$-Stirling numbers. Both correspond to the $\psi$
sequence choice in the $q$-Gauss form
$\langle\frac{1}{n_q!}\rangle_{n\geq 0}$. Here the
$q^\sim$-Stirling numbers of the second kind are introduced as
coefficients in the famous Newton interpolation formula (Liber
III, Lemma V, pp. 481-483 in [13]) now applied to the polynomial
sequence ${\langle e_n\rangle}_{n\geq0},\quad e_n(x)= x^n, n\geq
0$, i.e.

\begin{equation}\label{eq5}
x^n=\sum_{k=0}^{n}\Big\{{n \atop k}\Big\}^\sim_q \chi_{\underline
k}(x),\quad i.e. \quad \Big\{{n \atop k}\Big\}^\sim_q =
[0,1_q,2_q,...,k_q;e_n],
\end{equation}
where $\chi_{\underline k}(x)= x(x-1_q)(x-2_q)...(x-[k-1]_q)$, and
$$[x_0,x_1,x_2,...,x_k;f]=  \frac {
[x_0,x_1,x_2,...,x_{k-1};f]-[x_1,x_2,...,x_k;f]}{x_0-x_k}$$
denotes the $k-th$ divided difference with $$[x_0,x_1;f]= \frac {
f(x_0)-f(x_1)}{x_0-x_1}.$$
These two kinds of now classical
$q$-extensions of Stirling numbers of the second kind as defined
by (4) and (5)  are related in a simple way through re-scaling
[14]. They satisfy the known respective recurrences:

$$ \Big\{{{n+1}\atop k}\Big\}_q =
\sum_{l=0}^{n}\binom{n}{l}_q q^l\Big\{{l\atop {k-1}}\Big\}_q ;
n\geq 0, k\geq 1,$$

$$ \Big\{{{n+1}\atop k}\Big\}^\sim_q =
\sum_{l=0}^{n}\binom{n}{l}_q q^{l-k+1}\Big\{{l\atop
{k-1}}\Big\}^\sim_q;  n\geq 0, k\geq 1.$$ From the above it
follows immediately that corresponding $q$-extensions of  $B_n$
Bell numbers satisfy respective recurrences on their own:

$$ B_q(n+1) =
\sum_{l=0}^{n}\binom{n}{l}_q q^l B_q(l) ; n\geq 0,$$

$$ B^\sim_q(n+1) =
\sum_{l=0}^{n}\binom{n}{l}_q q^{l+1} \overline {B} ^\sim_q(l),
n\geq 0 $$ where
$$ \overline {B}^\sim_q(l)=
\sum_{k=0}^{l} q^k\Big\{{l\atop k}\Big\}^\sim_q.
$$
Different definitions via (4) and (5) equations  correspond
consequently to different $q$-counting  \cite{14}. With any other
choice out of countless choices of the $\psi$ sequence the
equation (5) becomes just the  definition of $\psi^\sim$-Stirling
(vide "Fibonomial-Stirling") numbers of the second kind $\Big\{{n
\atop k}\Big\}^\sim_\psi$  and then $\psi^\sim$-Bell numbers
$B^\sim_n(\psi)$ are defined as usual as sums of the corresponding
Stirling-like numbers  - where now $\chi_{\underline k}(x)$ in (5)
is to be replaced by $\psi_{\underline k}(x)=
x(x-1_\psi)(x-2_\psi)...(x-[k-1]_\psi)$. These
$\psi^\sim$-Stirling numbers of the second kind  for  $q$ case
identified as Comtet numbers in Wagner`s terminology \cite{15,14}
satisfy familiar recursion and are given by familiar formulas to
be presented soon. The extension of definition (4) of the $q$-
Stirling numbers of the second kind beyond this $q$-case i.e.
beyond the $\psi=\langle\frac{1}{n_q!}\rangle_{n\geq 0}$ choice is
not that mnemonic at all and the problem of recursion appears.
This part of alternative treatment is to be considered later on
after we exploit a little bit more some consequences of (4).
Namely - due to (4) one immediately notices that the expectation
value with respect to $q$-coherent state $ |\gamma> $ with $
|\gamma|=1 $ of the $n$-th power of the number operator is exactly
the popular $q$-Dobinski formula which can be given Blissard
calculus inspired umbral form - like in (3). It is enough to apply
to (4) $L_q$ - the average value functional for $q$-Poisson
distribution [1, 2]. The formula thus obtained may be also treated
as a definition of $q$-extended Bell numbers $ B_n(q)$

\begin{equation}\label{eq6}
         L_q(X_q^n)=B_n(q),\qquad n\geq 0
\end{equation}
due to the fact that this linear functional $L_q$ interpreted as
the average value functional for the random variable $X_q$ with
the $q$-Poisson distribution [1] ( $L_q(X_q )= 1$ ) satisfies

\begin{equation}\label{eq7}
L_q(X_q^{\underline{n}})=1, \qquad  n\geq 0.
\end{equation}
Then with the $q$-exponential polynomials
$$\varphi_n(x,q)=\sum_{k=0}^{n}\Big\{ {n \atop k}\Big\}_q x^k\quad\quad\quad\quad(q-exp-pol-I)$$
one obtains  for  $x=1$ in correspondence with  $L_q(X_q )= 1$ the
$q$ - formula of Dobinski type (compare with [12] and see (5.28)
in [15]):
$$ \varphi_n(1,q)=B_n(q)=L_q(X_q^n)=e_q^{-1}\sum_{k=0}^{\infty} \frac{k_q^n}{k_q!}
,\quad n\geq 0,\quad e_q^{-1}\equiv [exp_q\{1\}]^{-1}.$$ We arrive
to this simple conclusion using Jackson derivative difference
operator in place of  $D = d/dx$  in $q$ =1 case and the power
series generating function $G(t)$  for $q$-Poisson probability
distribution:

\begin{equation}\label{eq8}
p_n=[exp_q\lambda]^{-1}\frac{\lambda^n}{n_q!},
G(t)=\sum_{n\geq0}p_n t^n ,
\end{equation}
where $exp_q\lambda\equiv
exp_{\psi(q)}\lambda=\sum_{n\geq0}\psi_n(q)t^n,\quad
\psi_n(q)\equiv \frac{1}{n_q!}.$ Naturally

\begin{equation}\label{eq9}
p_n=[\frac{\partial_q^n G(t)}{n_q!}]_{t=0},\quad [\partial_q
G(t)]_{t=1} = 1\quad  for\quad   \lambda = 1.
\end{equation}
In order to arrive at the $q$-Dobinski formula apply (7) to (4)
with (6) in mind. (As for $\psi$-Poisson probability distribution
- see [1,2].) There are many ways leading to $q$-extended Stirling
numbers according to their weighted counting interpretation i.e.
various statistics are counted by $q$-Stirling numbers of the
second kind \cite{16}. For example $ w(\pi) = q^{cross(\pi)},
\quad w(\pi) = q^{inv(\pi)}$  from  \cite{17} give the
Carlitz-Gould-Milne $q$-Stirling numbers $\Big\{{n \atop
k}\Big\}_q $ - after being summed over  the set of $k$-block
partitions  while $w(\pi) = q^{nin(\pi)} $ from \cite{18} gives
rise to the Carlitz-Gould -Milne $q^\sim$-Stirling numbers $
\Big\{{n \atop k}\Big\}^\sim_q$ after being summed over the set of
$k$-block partitions (see also $maj$ and $maj^\sim$ statistics in
\cite{19} as well as other statistics [14] giving also both
extensions). The choice of weight $w(\pi) = q^{i(\pi)}$ \cite{20}
gives rise to another statistics counted by new kind of
$q$-extended Stirling numbers of the second kind. Let us consider
- for the sake of illustration this example from \cite{20} in more
detail. This is the example of weighted counting $\sum_{\sigma\in
\Pi_n}q^{i(\sigma)}$ of partitions of [n]. The weight  $w$  of
such a partition $\pi\in \Pi_n$ is determined by inversions` $i$
function in the form $w(\pi) = q^{i(\pi)}$. Here $\Pi_n$ denotes
the lattice of all partitions of the set [n] while $A_{n,k}$
denotes the family of all $k$ -block partitions. A $k$-block
partition $\pi\in A_{n,k}\subseteq\Pi_n$ is represented in the
standard form: $\pi = B_1/ B_2 /.../B_k$ with the convention that
$ max B_1< max B_2 <...,max B_k = n $. For $i\in[n]$  let $b_i$
denotes a number of a block to which $i$ pertains. Define an
inversion of partition $\pi$ to be  a pair $<i,j>$ such that
$b_i<b_j$ and $i>j$. The inversion set of $\pi$ is $I(\pi)=
\{<i,j> ;<i,j> \textit{is an inversion of}\quad \pi\}.$ Then
$i(\pi) = |I(\pi)|$ and the \textit{inversion} $q$-Bell numbers
are naturally defined as
$$ B^{inv}_n(q) = \sum_{\sigma\in \Pi_n}q^{i(\sigma)} =
\sum_{k \geq 0} \sum_{\pi\in A_{n,k}}{q^{i(\pi)}}$$ while
\textit{inversion} $q$-Stirling numbers of the second kind are
identified with
$$\Big\{{n\atop k}\Big\}^{inv}_q\equiv\sum_{\pi\in
A_{n,k}}{q^{i(\pi)}}.$$ The \textit{inversion} $q$-Bell number
$B^{inv}_n(q)$ is the generating function for $I(s) =$ the number
of all partitions from $\Pi_n$ with $s$ of the above inversions
because
$$B^{inv}_n(q) = \sum_{\sigma\in \Pi_n}q^{i(\sigma)} =
\sum_{s\geq 0} q^s \sum_{\pi\in \Pi_n , i(\pi)=s }1\equiv
\sum_{s\geq 0} I(s)q^s .$$ Recursions for both \textit{inversion}
$q$-Bell numbers and \textit{inversion} $q$-Stirling  numbers of
the second kind are not difficult to be derived. Also in a natural
way the \textit{inversion} $q$-Stirling numbers of the second kind
from [20] satisfy a $q$-analogue of the standard recursion for
Stirling numbers of the second kind to be written via mnemonic
adding "$q$" subscript to the binomial and second kind Stirling
symbols in the the standard recursion formula i.e.

$$ \Big\{{{n+1}\atop k}\Big\}^{inv}_q =
\sum_{l=0}^{n}\binom{n}{l}_q \Big\{{{n-l}\atop
{k-1}}\Big\}^{inv}_q  ;\quad n\geq 0 , k\geq 1.$$ Another
$q$-extended Stirling numbers much different from Carlitz
"$q$-ones" were introduced in the reference \cite{21} from where
one infers \cite{22} the \textit{cigl}-analog of (6). Let $\Pi$
denotes the lattice of all partitions of the set $\{0,1,..,n-1\}$.
Let $\pi\in \Pi $  be represented by blocks   $ \pi =\{B_o ,B_1
,...B_i ,...\}$, where $B_o$ is the block containing zero:  $0\in
B_o$. The weight adapted by Cigler defines weighted partitions`
counting according to the content of $B_o$. Namely

$$w(\pi)= q^{cigl(\pi)}, cigl(\pi)=\sum_{l\in {B_0}}l ,
\sum_{\pi\in A_{n,k}}{q^{cigl(\pi)}}\equiv\Big\{{n\atop
k}\Big\}^{cigl}_q .$$ Therefore $ \sum_{\pi\in\underline{}
\Pi}{q^{cigl(\pi)}}\equiv{B_n(q)}$. Here $A_{n,k}$ stays for
subfamily of all $k$-block partitions. With the above relations
one has defined the \textit{cigl}-$q$-Stirling and the
\textit{cigl}-$q$-Bell numbers. The \textit{cigl}-$q$-Stirling
numbers of the second kind are expressed in terms of $q$-binomial
coefficients and $q =1$ Stirling numbers of the second kind
\cite{21} as follows

$$ \Big\{{{n+1}\atop k}\Big\}^{cigl}_q =
\sum_{l=0}^{n}\binom{n}{l}_q q^{\binom{n-l+1}{2}}\Big\{{{n-l}\atop
{k-1}}\Big\}^{cigl}_q ; n\geq 0 , k\geq 1.$$ As seen above these
are \textbf{new} $q$-extended Stirling numbers. The corresponding
\textit{cigl}-$q$-Bell numbers recently have been equivalently
defined via \textit{cigl}-$q$-Dobinski formula \cite{22} - which
now in more adequate notation reads:

$$L(X^{\overline{q^n}})=\overline{B}_n(q),\qquad n\geq 0,
X^{\overline{q^n}}\equiv X(X+q-1)...(X-1+q^{n-1}).$$ The above
\textit{cigl}-$q$-Dobinski formula is interpreted as the average
of this specific  $n-th$ \textit{cigl}-$q$-power random variable
$X^{\overline{q^n}}$ with the $q = 1$ Poisson distribution such
that $ L(X)=1.$ For that to see use the identity by Cigler
\cite{21}
$$
x(x-1+q)...(x-1+q^{n-1})=\sum_{k=0}^{n}\Big\{ {n \atop
k}\Big\}^{cigl}_q  x^{\underline k}.
$$

\subsection*{2. Beyond the $q$-extensions realm}

The further consecutive $\psi$-umbral extension of Carlitz-Gould
$q$-Stirling numbers $\Big\{{n\atop k}\Big\}_q $  and
$\Big\{{n\atop k}\Big\}^\sim_q $ is realized two-fold way - one
of which leads to a surprise  in  contrary to the other way. \\

\subsection*{2.1. The first way}
\textbf{The first} "easy way" consists in almost mnemonic
sometimes replacement of $q$ subscript by $\psi$ after having
realized that in equation (5) we are dealing  with the specific
case of the so called Comtet numbers \cite{14,15} (Comtet L. in
\emph{Nombres de Stirling generaux et fonctions symtriques}  C.R.
Acad. Sci. Paris, Series A, 275 (1972):747-750 formula (2) refers
to \textbf{Wronski}). This array of Stirling-like numbers
$\Big\{{n \atop k}\Big\}^\sim_{\psi}$ - \textit{"alephs de
Wronski"} as Comtet refers to it  or these Comtet numbers in
terminology of Wagner [14, 15] or as a matter of fact [13] these
\textbf{Newton} interpolation coefficients  for  $e_n, n\geq0$
i.e. divided differences $[0,1_{\psi},2_{\psi},...,k_{\psi};e_n]$
are defined accordingly as such coefficients - below.

\begin{equation}\label{eq10}
x^n=\sum_{k=0}^{n}\Big\{{n \atop k}\Big\}^\sim_{\psi}
\psi_{\underline k}(x),\quad n\geq0,
\end{equation}
i.e. equivalently  (recall that  $e_n(x)= x^n, n\geq 0$)

$$\quad \Big\{{n \atop k}\Big\}^\sim_{\psi} =
[0,1_{\psi},2_{\psi},...,k_{\psi};e_n]=\sum_{l=0}^{k}\frac
{e_n(l_{\psi})}{{\psi}^`_{\underline{k+1}}(l_{\psi})},\quad
n\geq0,\quad\quad \quad(Newton)$$ where
$$\psi_{\underline k}(x)=
x(x-1_{\psi})(x-2_{\psi})...(x-[k-1]_{\psi})$$ and
$\psi^`_{\underline s}$ denotes the first derivative. Let then $f
= \langle f_n \rangle_{n\geq0}$ be an arbitrary sequence of
polynomials. In the following we shall call $S(f;n,k)$ defined
below
$$[d_0,d_1,d_2,...,d_k;f_n] \equiv S(f;\langle d_l \rangle_{l\geq0},n,k)\quad
\quad (N-W-C\quad Stirling)$$ the Newton-Wronski-Comptet
Stirling numbers (N-W-C  for short)- compare with Appendix A.2.\\
The $\psi^\sim $ - Stirling numbers  $\Big\{{n \atop
k}\Big\}^\sim_{\psi}$ defined by  (10) are specification of
$N-W-C$ Stirling array for which we naturally define $\psi^\sim
$-exponential polynomials $\varphi_n(x,\psi)$ as follows
$$\varphi_n^\sim (x,\psi)
=\sum_{k=0}^{n}[0,1_{\psi},2_{\psi},...,k_{\psi};e_n]x^k, \quad
n\geq0.\quad\quad(\psi^\sim-exp-pol)$$ \textbf{Note }the trivial
but important fact that in the N-W-C Stirling numbers case we are
dealing with not equidistant nodes` interpolation in general and
note that ($Rescal$) from the subsection 2.2. below  is no more
valid beyond $q$-extension case - both with an impact on the way
to find out the Dobinski-like formulae - see more below.\\
As a consequence of (10) we have "for granted" the following
extensions of recurrences for Stirling numbers of the second kind:

\begin{equation}\label{eq11}
\Big\{{{n+1}\atop k}\Big\}^\sim_{\psi} = \Big\{{n\atop
{k-1}}\Big\}^\sim_{\psi} + k_{\psi}\Big\{{n\atop
k}\Big\}^\sim_{\psi} ;\quad n\geq 0, k\geq 1,
\end{equation}
where \quad $\Big\{{n\atop 0}\Big\}^\sim_{\psi}=
\delta_{n,0},\quad\Big\{{n\atop k}\Big\}^\sim_{\psi}=0 ,\quad k>n
;\quad $\quad and the recurrence for ordinary generating function
reads

\vspace{1mm}

\begin{equation}\label{eq12}
G^\sim_{k_{\psi}}(x)=\frac{x}{1-k_{\psi}}G^\sim_{k_{\psi}-1}(x) ,
\quad k\geq 1,
\end{equation}
where naturally
$$G^\sim_{k_{\psi}}(x)= \sum_{n\geq 0}\Big\{{n\atop
k}\Big\}^\sim_{\psi}x^n,\quad k\geq 1 $$ from where one infers
that

\begin{equation}\label{eq13}
G^\sim_{k_{\psi}}(x)=\frac{x^k}{(1-1_{\psi}x)(1-2_{\psi}x)...(1-k_{\psi}x)}
\quad, \quad k\geq 0.
\end{equation}
Hence we arrive in the standard extended  text-book way \cite{22}
at the following explicit \textbf{new} formula (compare with (2.3)
in [15])

\begin{equation}\label{eq14}
[0,1_{\psi},2_{\psi},...,k_{\psi};e_n]=\Big\{{n\atop
k}\Big\}^\sim_{\psi} = \frac
{1}{k_{\psi}!}\sum_{r=1}^{k}(-1)^{k-r}\binom{k_{\psi}}{r_{\psi}}r_{\psi}^n
; \quad n\geq k \geq 0,
\end{equation}
where
$$\sum_{r=1}^{k}(-1)^{k-r}\binom{k_{\psi}}{r_{\psi}}r_{\psi}^n ;
\quad n,k\geq 0 $$ is readily recognized as the $\psi$-extension
of the formula for surjections in its - after inclusion-exclusion
principle had been applied - form.

\vspace{1mm} \noindent Expanding the right hand side of (13)
results in another explicit formula for these $\psi$-case
Newton-Wronski-Comtet array of Stirling numbers of the second kind
i.e. we have

\begin{equation}\label{eq15}
\Big\{{n\atop k}\Big\}^\sim_{\psi} = \sum_{1\leq i_1 \leq
i_2\leq...\leq i_{n-k}\leq
k}(i_1)_{\psi}(i_2)_{\psi}...(i_{n-k})_{\psi}; \quad \quad n \geq
k\geq 0
\end{equation}
or equivalently  (compare with [13, 14])

\begin{equation}\label{eq16}
\Big\{{n\atop k}\Big\}^\sim_{\psi}=\sum_{d_1+ d_2+...+d_k =
n-k,\quad d_i\geq
0}1_{\psi}^{d_1}2_{\psi}^{d_2}...k_{\psi}^{d_k};\quad n \geq k\geq
0.
\end{equation}
N-W-C case $\psi^\sim$ - \textbf{Stirling numbers} of the second
kind being defined equivalently by (10), (\textit{Newton}), (14),
(15) or (16) yield N-W-C case $\psi^\sim$ - \textbf{Bell numbers}

$$ B^\sim_n(\psi)=\sum_{k=0}^{n} \Big\{{n\atop
k}\Big\}^\sim_{\psi}=\sum_{k=0}^{n}[0,1_{\psi},2_{\psi},...,k_{\psi};e_n]
,\qquad n\geq 0 \quad\quad\quad\quad(B^\sim).$$ Naturally
$\exists$! functional $L^\sim $ such that on the basis of
persistent root polynomials $\psi_{\underline k}(x)$ it takes the
value $1$:

$$L^\sim (\psi_{\underline k}(x))=1,\quad k\geq 0. $$
Then from (10) we get an analog of (3)

$$B^\sim_n (\psi)= L^\sim (x^n) \quad\quad\quad\quad\quad\quad\quad\quad\quad\quad\quad(L^\sim). $$

\noindent \textit{Problem}: which distribution the functional
$L^\sim $ is related to is an open technical question by now. More
- the recurrence  for  $B^\sim_n(\psi)$ is already quite involved
and complicated for the $q$-extension case (see: the first
section)- and no acceptable readable form of recurrence  for the
$\psi$-extension case is known to us by now.\\
\textbf{Nevertheless} after adapting  the standard text-book
method \cite{23} we have the following  formulae for two variable
ordinary generating function for  $\Big\{{n\atop
k}\Big\}^\sim_{\psi}$ Stirling numbers of the second kind and the
$\psi$-exponential generating function for $ B^\sim_n(\psi)$  Bell
numbers
\begin{equation}\label{eq17}
C^\sim_{\psi}(x,y) = \sum_{n\geq 0} \varphi^\sim_n (\psi,y)x^n,
\end{equation}
where the $\psi$- exponential polynomials $ \varphi^\sim_n
(\psi,y)$ $$ \varphi^\sim_n (\psi,y)=\sum_{k=0}^{n} \Big\{{n\atop
k}\Big\}^\sim_{\psi}y^k$$ do satisfy the recurrence (compare with
formulas (28) in Touchard`s [24] from 1956)

$$ \varphi^\sim_n (\psi,y)=  [y(1+\partial_{\psi}]\varphi^\sim_{n-1}(\psi,y) \qquad n\geq 1 ,$$
hence
$$ \varphi^\sim_n (\psi,y)=  [y(1+\partial_{\psi}]^n 1,\quad \qquad n\geq 0.$$
The linear operator $\partial_{\psi}$ acting on the algebra of
formal power series is being called (see: [1, 2] and references
therein) the "$\psi$-derivative" as $\partial_{\psi} y^n =
n_{\psi}y^{n-1}.$\\
The $\psi^{\sim}$ - exponential generating function

$$
B^{\sim}_{\psi}(x)= \sum_{n\geq
0}B^\sim_n(\psi)\frac{x^n}{n_{\psi}!}\quad\quad\quad\quad\quad(\psi^{\sim}-
e.g.f.)
$$
for $ B^{\sim}_n(\psi)$ Bell numbers - after cautious adaptation
of the method from the Wilf`s generatingfunctionology book  [23]
can be seen to be given by the following \textbf{new }formula

\begin{equation}\label{eq18}
B^{\sim}_{\psi}(x)= \sum_{r\geq 0}\epsilon(\psi,r)
\frac{e_{\psi}[r_{\psi}x]}{r_{\psi}!}
\end{equation}
where (see: [1,2] and references therein)
$$e_{\psi}(x) =
\sum_{n\geq 0}\frac{x^n}{n_{\psi}!}$$ while
\begin{equation}\label{eq19}
\epsilon(\psi,r)=\sum_{k=r}^{\infty} \frac{(-1)^{k-r}}{(k_{\psi}-
r_{\psi})!}
\end{equation}
and the \textbf{new } Dobinski - like formula for the
$\psi$-extensions here now reads

\begin{equation}\label{eq20}
B^{\sim}_n (\psi)= \sum_{r\geq 0}\epsilon(\psi,r)
\frac{r_{\psi}^n}{r_{\psi}!}.
\end{equation}
The $\psi^\sim$-exponential polynomials are therefore given
correspondingly by
$$ \varphi^\sim_n (\psi,x)=\sum_{r\geq 0}\epsilon(\psi,r)
\frac{r_{\psi}^n}{r_{\psi}!}x^r.\quad\quad\quad\quad\quad(\psi^\sim-exp-pol-II)$$
In the case of Gauss $q$-extended  choice of
$\langle\frac{1}{n_q!}\rangle_{n\geq 0}$ admissible sequence of
extended umbral operator calculus equations (19) and (20) take the
form

\begin{equation}\label{eq21}
\epsilon(q,r)=\sum_{k=r}^{\infty}
\frac{(-1)^{k-r}}{(k-r)_q!}q^{-\binom {r}{2}}
\end{equation}
and the\textbf{ new} N-W-C case $q^\sim$-Dobinski  formula is
given by

\begin{equation}\label{eq22}
B^{\sim}_n (q)= \sum_{r\geq 0}\epsilon(q,r) \frac{r_q^n}{r_q!},
\end{equation}
which for $ q=1$ becomes the Dobinski formula from 1887 [4]. Note
the appearance of re-scaling factor  $q^{-\binom {r}{2}}$ in (21).
In its absence we would get \textbf{not} $q^\sim$-Dobinski  but
$q$-Dobinski formula

$$B_n (q)= \frac {1}{\exp_q(1)}\sum_{0\leq
k}\frac{k_q^n }{k_q!}\quad\quad\quad\quad\quad( q-Dobinski)$$ -
see [15] and formula (5.28) there coinciding with   N-W-C  case of
Dobinski formula after re-scaling  in correspondence with (Rescal)
below in subsection 2.2. Correspondingly we would the have
\textbf{ not} $(q^\sim-exp-pol)$ formula but $(q-exp-pol)$
formula:

$$\varphi_n(x,q)=\sum_{k=0}^{n}q^{\binom{k}{2}}[0,1_q,2_q,...,k_q;e_n]x^k=
\sum_{k=0}^{n}\Big\{ {n \atop k}\Big\}_q x^k. \quad\quad(q-exp-pol)$$\\
\textbf{The interpretation problem.} Combinatorial interpretations
of the known up to now various $q$-extensions of Stirling numbers
of both kinds - are briefly reported on in  the Appendix. The
problem of how eventually one might interpret - beyond the
q-extensions` realm -  for example the $\psi^\sim$-Dobinski
formulae (20) and (22) also in the Rota-like way represented  here
by equation (3)  we leave \textbf{opened} - see the discussion in
Appendix A.2 and A.3.II. Naturally there exist a unique linear
functional $L^\sim_{\psi}$ such that
$$L^\sim_{\psi}(\psi_{\underline
k}(x)= x(x-1_q)(x-2_q)...(x-[k-1]_{\psi}))=1,\quad k\geq 0.$$ It
is also to be noted that in the exceptional  case of
$q$-extensions and \textbf{only for $q$-extension} we have
equivalence of $\Big\{{r\atop l}\Big\}^\sim_{\psi}$ and
$\Big\{{r\atop l}\Big\}_{\psi}$ by re-scaling. \\
For the latter $q$-Stirling numbers we have Dobinski formula and
simultaneously $q$-Poisson average functional interpretation as
represented by the definition (6). Namely - recall  the fact that
the linear functional $L_q$ there is interpreted as the average
value functional for the random variable  $X_q$ with the
$q$-Poisson distribution  which is specific case of the
$\psi$-Poisson distributions from [2, 1]. Note again that this
re-scaling takes place for $q$-extensions and \textbf{only for
$q$-extension}. This is so because the relation

$$n_{\psi}- k_{\psi} = f(k)(n-k)_{\psi}, \quad 1_{\psi}=1$$
holds for and only for $f(1)\equiv q$ when it becomes the identity
(2.6) from [25] i.e. $f(k)=q^k$. It is our conviction that this is
the very reason that $q$-extensions seem to  appear as a kind of
"a bifurcation point" in the domain of $\psi$-umbral extensions.
This \textbf{conviction} is supported by the corresponding
considerations in [26]- section 3 - on possibility of
${\psi}$-analogue of the so called "quantum" $q$-plane formulation
of $q$-umbral calculus. \vspace{1mm} The parallel treatment  of
the Newton-Wronski-Comtet
$\left[ {{\begin{array}{*{20}c} {n}\\
{k}\end{array}} } \right]^\sim_\psi$ Stirling numbers of the
\emph{first kind }is now not difficult (consult [25, 12, 8] for
example and Wagner`s recent treatment of the well established
$q$-case in [14]).

\vspace{1mm}

\noindent In the inversion-dual way to our equation (10) above we
define the $\psi^\sim$-Stirling numbers of the first kind as
coefficients in the following expansion

\begin{equation}\label{eq23}
\psi_{\underline k}(x)=\sum_{r=0}^{k}\left[ {{\begin{array}{*{20}c} {k}\\
{r}\end{array}} } \right]^\sim_\psi x^r
\end{equation}
where - recall $\psi_{\underline k}(x)=
x(x-1_{\psi})(x-2_{\psi})...(x-[k-1]_{\psi}).$ (Attention: see
equations (10)-(16) in [8] and note the difference with the
present definition). Therefore from the above we  infer that

\begin{equation}\label{eq24}
\sum_{r=0}^{k}\left[ {{\begin{array}{*{20}c} {k}\\
{r}\end{array}} } \right]^\sim_\psi \Big\{{r\atop
l}\Big\}^\sim_{\psi}= \delta_{k,l}.
\end{equation}
Another natural counterpart to $\psi^\sim$-Stirling numbers of the
second are $\psi^c$- Stirling numbers  of the first kind defined
here down as coefficients in the following expansion ("\textit{c}"
because of cycles in non-extended case)
\begin{equation}\label{eq25}
\psi_{\overline k}(x)=\sum_{r=0}^{k}\left[ {{\begin{array}{*{20}c} {k}\\
{r}\end{array}} } \right]^c_\psi x^r
\end{equation}
where - now  $\psi_{\overline k}(x)=
x(x+1_{\psi})(x+2_{\psi})...(x+[k-1]_{\psi}).$ These are to be
studied elsewhere.\\
\textbf{On interpretation.} For possible \textbf{unified}
combinatorial interpretations of binomial coefficients of both
kinds, the Stirling numbers of both kinds and the Gaussian
coefficients of the first and second kind - i.e for the specific
choices of $\psi = \langle \frac{1}{n_{\psi!}}\rangle_{n\geq 0}$ -
see [27, 28]. As for q-analogue of Stirling cycle numbers see [29]
and  Sect. 5.3. in [30]. The problem of eventual combinatorial
interpretation of other  $\psi$-extensions (vide Fibonomial - for
example) - remains opened.\\

\subsection*{2.2. The second way.}
We shall come over now to inspect the outcomes of \textbf{the
second way} - announced at the start of this section - after
having realized that in the equation (4) we are \textbf{not}
dealing straightforwardly  with Newton-Wronski-Comtet array of
Stirling-like numbers \cite{13, 14, 15}  - except for the $q=1$
case "extension"- of course. Though it is to be noted that still
the re-scaling takes place for $q$-extensions and \textbf{only for
$q$-extension}. Not beyond. (This is so because the relation
$$n_{\psi}- k_{\psi} = f(k)(n-k)_{\psi}, \quad 1_{\psi}=1$$ holds
for and only for $f(1)\equiv q$ when it becomes the identity (2.6)
from Gould`s [25] i.e. $f(k)=q^k$). Thus after the above  Gould
re-scaling we would recover  Newton-Wronski-Comtet array of
Stirling-like numbers re-scaled - anyhow i.e.
$$y^n=\sum_{k=0}^{n}q^{-\binom{k}{2}} \Big\{{n \atop k}\Big\}_q y^{\underline k},\quad\quad\quad\quad(Rescal)$$
where $ y=x_q=\frac{1-q^x}{1-q}$ and $\psi_{\underline k}(y)=
y(y-1_q)(y-2_q)...(y-[k-1]_q).$ \vspace{1mm} At first let us
recall that the definition (4) of $q$-Stirling numbers of the
second kind $\Big\{{n \atop k}\Big\}_q$ is equivalent to the
definition by the recursion

\begin{equation}\label{eq26}
\Big\{{{n+1}\atop k}\Big\}_q = q^{k-1}\Big\{{n\atop {k-1}}\Big\}_q
+ k_q\Big\{{n\atop k}\Big\}_q ;\quad n\geq 0 , k\geq 1,
\end{equation}
where \quad $\Big\{{n\atop 0}\Big\}_q =
\delta_{n,0},\quad\Big\{{n\atop k}\Big\}_q=0 ,\quad k>n.$

\vspace{1mm}

\noindent These in turn is equivalent to (just use the standard
$Q$-Leibniz rule [1, 2, 31] for Jackson derivative $\partial_q$)

\begin{equation}\label{eq27}
(\hat{x} \partial_q)^n =\sum_{k=0}^{n}\Big\{{n \atop k}\Big\}_q
\hat{x}^k \partial_q^k,
\end{equation}
where \quad $\Big\{{n\atop 0}\Big\}_q =
\delta_{n,0},\quad\Big\{{n\atop k}\Big\}_q=0 ,\quad k>n.$ Here
$\hat{x}$ denotes the multiplication by the argument of a
function. The formula (27) is a special case of the typical for
GHW algebra [32, 1, 2, 33] expression investigated by Carlitz in
1932 [34] (compare with GHW formulae (1), (31), (32) in [7] and
see also [35]). \vspace{1mm} The idea now is to extend eventually
the definition by equation (27)  \textit{via} \textbf{replacing}
$q$-extended operators by the corresponding $\psi$-extended
elements of the Graves-Heisenberg-Weyl (\textit{GHW}) algebra
representation [32, 1, 2, 33]. However we note at once (see:
Appendix for particulars of the up-side down notation) that the
two following observations hold.\\
\textbf{Observation 2.1} The equivalent definitions  (28) and (29)

\begin{equation}\label{eq28}
( \hat{x} \partial_\psi)^n =\sum_{k=0}^{n}\Big\{{n \atop
k}\Big\}_\psi \hat{x}^k
\partial_\psi^k
\end{equation}
where \quad $\Big\{{n\atop 0}\Big\}_\psi =
\delta_{n,0},\quad\Big\{{n\atop k}\Big\}_\psi=0 ,\quad k>n$ and

\begin{equation}\label{eq29}
x_\psi^n=\sum_{k=0}^{n}\Big\{{n \atop k}\Big\}_\psi
x_\psi^{\underline k}
\end{equation}
lead to the one first order recurrence of the type (26)
\textbf{for and only for $q$-extension}.\\
This again is so because the relation
$$n_{\psi}- k_{\psi} = f(k)(n-k)_{\psi} , \quad 1_{\psi}=1$$
holds for and only for  $f(1)$ i.e. for  {$q$-extension}, where it
becomes the identity (2.6) from [25] i.e. $f(k)=q^k$. \vspace{1mm}
The next observation now comes as a would be surprise.\\
\textbf {Observation 2.2} The equivalent definitions (28) and (29)
have no non-trivial realizations   beyond the $q$-extension case.

\vspace{1mm}

In order to arrive at this observation let us act appropriately on
$x^N$\quad $N\geq 0$\quad monomials by both sides of the
\textit{GHW} algebra representation definition (28) thus getting
an \textit{infinite} \textbf{sequence }\textit{of recurrences}

\begin{equation}\label{eq30}
\sum_{k\geq0}\Big\{{{n+1} \atop k}\Big\}_\psi N_\psi^{\underline
k}=\sum_{k\geq0}\Big\{{n \atop k}\Big\}_\psi N_\psi
N_\psi^{\underline k},\quad N\geq 0,
\end{equation}
\textit{with no nontrivial solutions} as spectacularly evident
with the choice of - for example - Fibonomialy-extended sequence
$\langle\frac{1}{F_n!}\rangle_{n\geq 0}$ ($\langle F_n \rangle$ -
Fibonacci sequence ) \textbf{unless}
$\psi=\langle\frac{1}{n_q!}\rangle_{n\geq 0}$. And for this and
only for this choice  $\psi=\langle\frac{1}{n_q!}\rangle_{n\geq
0}$ we have

$$N_q=q^k((N-k)_q + k_q$$
which after being applied in (30) results in \textbf{one
}recurrence which is exactly  the recurrence (26).

\vspace{1mm }

As expected - the equation (29) becomes equivalent to the one
first order recurrence of the type (26) \textbf{for and only for
$q$-extension}.\\
\textbf{Closing remark.} We see that the Carlitz-Gould
$q$-Stirling numbers $\Big\{{n\atop k}\Big\}_q $ make $q$-umbral
extension to appear as a kind of "a bifurcation point"  in the
domain of respective $\psi$-umbral extensions. This in statu
nascendi conviction is also supported by the corresponding
considerations in [26]- section 3 - considerations about
possibility of ${\psi}$-analogue of the so called "quantum"
$q$-plane formulation of $q$-umbral calculus. As for eventual
second way`s $\psi$-extensions \emph{beyond} the $q$-extension
case where the $rescaling$ does not take place we are left with an
opened problem how to eventually find the way to get round this
inspiring obstacle. The selective comparison of the presented
umbral extensions of Stirling numbers, Bell numbers and
Dobinski-like formulas with other existing extensions (as well as
relevant information in brief)  serves the purpose of seeking
analogies and is to be find in the Appendix that follows now.

\subsection*{Appendix - for remarks, discussion and brief comparative review of ideas.}
\subsection*{A.1. Notation.}
The necessary commutation relations` representation for the $GHW$
(Graves-Heisenberg-Weyl) algebra generators is provided in [31, 1,
2, 33]. Applications of these might be worthy of the further study
[25, 36, 37]. The simplicity of the first steps to be done while
identifying general properties of such $\psi$-extensions consists
in notation i.e. here - in writing objects of these extensions in
mnemonic convenient \textbf{upside down notation} \cite{1},
\cite{2}

\begin{equation}\label{eq31}
\frac {\psi_{(n-1)}}{\psi_n}\equiv n_\psi,\quad
n_\psi!=n_\psi(n-1)_\psi!, \quad n>0, \quad x_{\psi}\equiv \frac
{\psi{(x-1)}}{\psi(x)},
\end{equation}

\begin{equation}\label{eq32}
x_{\psi}^{\underline{k}}=x_{\psi}(x-1)_\psi(x-2)_{\psi}...(x-k+1)_{\psi}
\end{equation}

\begin{equation}\label{eq33}
x_{\psi}(x-1)_{\psi}...(x-k+1)_{\psi}=
\frac{\psi(x-1)\psi(x-2)...\psi(x-k)} {\psi(x)
\psi(x-1)...\psi(x-k +1)}.
\end{equation}
If one writes the above in the form $x_{\psi} \equiv \frac
{\psi{(x-1)}}{\psi(x)}\equiv \Phi(x)\equiv\Phi_x\equiv x_{\Phi}$,
one sees that the name upside down notation is legitimate. You may
consult [1, 2, 26, 33, 36, 37] for further development and
usefulness of this notation. In this notation the $\psi$-extension
of binomial incidence coefficients read familiar:
$$\Big({n \atop k}\Big)_{\psi}=
\frac{n_{\psi}!}{k_{\psi}!(n-k)_{\psi}!}=\Big({n \atop {n-k}}\Big
)_{\psi}.$$

\subsection*{A.2. Discussion, remarks, questions.}
$q$-umbral extensions are expected to be of distinguished
character - also due to what was stated in Section 2. Being so
they pay to us with simplicity of formulae and elegance of $q$ -
weighting combinatorial interpretations and thus various
statistics of the combinatorial origin.

\vspace{1mm}\noindent Because of that and because of the
importance of $q$-umbral extensions in coherent states mathematics
we here
adjoin a remark on simplicity based on the remark of Professor Cigler (in private).\\
Namely  Katriel indicates  in the very important source  paper
\cite{6} that his derivation of the Dobinski formula is the
simplest. And really it is simple and wise. Possibly then this may
be occasionally and profitably confronted with  the also extremely
simple derivation by Cigler (see p. $104$ in [38]) based on
GHW-algebra properties. \vspace{1mm} Let then $\hat{x}$ denotes
the multiplication by $x$ operator while $D$ denotes
differentiation - both acting on the prehilbert space $P$ of
polynomials. Then due to the recursion for Stirling numbers of the
second kind and the identity (operators act on $P$)

$$ \hat{x}(D+1)\equiv \frac {1}{\exp{(x)}}(\hat{x}D)\exp{(x)}$$
one defines in GHW - algebra manner the \textbf{exponential
polynomials}

$$n\geq0,\quad\quad\quad\varphi_n(x)=\sum_{k=0}^{n}\Big\{ {n \atop k}\Big\}x^k \quad
\quad\quad\quad\quad\quad\quad(ExPol)$$ introduced by Acturialist
J.F. Steffensen [39, 40] (see: Bell`s "Exponential polynomials" in
umbra-symbolic language [41] p. 265 and his symbolic formula (4.7)
for now Bell numbers). These exponential polynomials were
substantially investigated by Touchard in Blissard umbra-symbolic
language [24]. Here now comes the GHW-definition [38] of these
basic polynomials

$$ \varphi_n(x)= \frac {1}{\exp{(x)}}(\hat{x}D)^n \exp{(x)}$$
resulting in the formula which becomes Dobinski one for  $x=1$
i.e.
$$\varphi_n(x)= \frac {1}{\exp{(x)}}\sum_{0\leq k}\frac{k^n
 x^k}{k!}.$$
\textbf{Note:} The $q$-case as well as $\psi$-case \textbf{formal}
mnemonic counterpart formulae are automatically arrived at with
the mnemonic attaching of $q$ or $\psi$ indices to nonnegative
numbers \cite{1,2} - vide:
\begin{equation}\label{eq34}
\varphi_n(x,\psi)= \frac {1}{\exp_{\psi}{(x)}}\sum_{0\leq
k}\frac{k_{\psi}^n x^k}{k_{\psi}!}
\end{equation}
which for  $\psi = \langle\frac{1}{n_q!}\rangle_{n\geq 0}$ and
$x=1$ becomes the well known $q$-Dobinski formula as of course
\quad $\varphi_n(x=1,q)= B_n(q)$ - see in [15] the formula (5.28)
and note that this is \textbf{not} $q^\sim$-Dobinski formula (22)
as noticed right after (22). As for eventual second way`s
$\psi$-extensions beyond the $q$-extension case where the
$rescaling$ does not take place - we are left with an opened
problem how to eventually find the way to get round this inspiring
obstacle. Perhaps instead of  the second beyond the $q$-extension
way we might follow Alexander the great in his Gordian Knot
problem solution  and define $S(\psi,n,k)$ as follows (whenever
one may prove that the object being defined is really a
polynomial):

$$\varphi_n(x,\psi)=\sum_{k=0}^{n}S(\psi,n,k)x^k = \frac{1}{\exp_{\psi}(x)}\sum_{0\leq
k}\frac{k_{\psi}^n x^k}{k_{\psi}!}. \quad\quad(S(\psi)-exp-pol)$$
An alternative  good idea perhaps would be an attempt to
$\psi$-extend the celebrated Newton interpolation formula ( use
$\partial_{\psi}$ instead $D$, then $exp_{\psi}$ instead of $exp$
and then you will be faced with $\psi$-Leibniz rule application
problem though... see  [1, 2, 33] for Leibnitz rules). Let us then
make - also for the sake of comparison with existing knowledge -
let us then make us wonder on the intrinsic presence and
assistance of Newton interpolation  which corresponds to the first
"easy"
way as described in Subsection 2.1.\\
\textbf{The intrinsic presence and assistance of Newton
interpolation} formula in derivation of Dobinski formula for
exponential polynomials and their binomial analogues was
underlined and used in [42] for specific presentation of the $q=1$
case from the umbral point of view of the classical finite
operator calculus. In [42] a Dobinski-like formula was derived
being as a matter of fact the particular ("binomial") case of
formula (30) from Touchard`s 1956 year paper [24]. In more detail.
Choosing any binomial polynomial sequence ${\langle
b_n\rangle}_{n\geq0}$ consider its Newton interpolation formula
$$b_n(x) =\sum_{k=0}^{n}[0,1,2,...,k;b_n]x^{\underline k}. $$
Then apply an umbral operator sending the binomial basis ${\langle
x^{\underline n}\rangle}_{n\geq0}$  of delta operator $\Delta$ to
the binomial basis $\langle x^n\rangle_{n\geq0}$ of delta
operator $D$. Then use\\
$$[0,1,2,...,k;b_n]= \frac{\Delta^k b_n|_{x=0}}{k!}=
\sum_{l=0}^{k}\frac
{(-1)^{k-l}b_n(l)}{(k-l)!l!}\quad(Newton-Stirling)$$ so as to
arrive (thanks to binomial convolution)  at Dobinski like formula
from [42]  i.e.

$$
b_n(\varphi(x)) = \frac
{1}{\exp(x)}\sum_{k=0}^{\infty}\frac{b_n(k)
 x^k}{k!},
 $$
where $\varphi$ is the umbral symbol satisfying [24]
$$\varphi_{n+1}= x(\varphi + 1)^n,\quad \varphi^{\underline k}= x^k. \quad\quad\quad\quad\quad (Touchard)$$
In order to see that this is just the particular ("binomial") case
of umbra-symbolic formula (30) from Touchard`s 1956 year paper
[24] just choose in Touchard formula (30) the \textit{arbitrary
}polynomial $f$ to be any binomial one  $b_n = f $. Then
$f(\varphi)=b_n(\varphi)=b_n(\varphi(x))$ is binomial also and we
have
$$f(\varphi) = \frac {1}{\exp(x)}\sum_{k=0}^{\infty}\frac{b_n(k)
x^k}{k!}.\quad\quad\quad\quad\quad\quad(Dobinski-Touchard)$$
\textbf{Equidistant nodes} Newton`s interpolation array of
coefficients $[0,1,2,...,k;b_n]$ - here the connection constants
of the general exponential polynomial $p_n(x) = b_n(\varphi(x))$
are to be called in the following the \textbf{Newton-Stirling}
numbers of the second kind and are consequently given by
$$ p_n(x) = \frac {1}{\exp(x)}\sum_{k=0}^{\infty}\frac{b_n(k)
x^k}{k!}= \sum_{k=0}^{n}[0,1,2,...,k;b_n]x^k
,\quad\quad\quad\quad(N-S-Dob)$$ where $\langle b_n
\rangle_{n\geq0}$ is any sequence of polynomials.  These are - in their turn - the special case of N-W-C Stirling numbers. \\
\textbf{Coherent States` Example I.}  Take the $b_m(x)=f(x)$  in
the (\textit{Dobinski-Touchard}) formula to be of the form
resulting from normal ordering problem (see A.3.II. - below) i.e.
let (see: [10])
$$f(x)= b_{ns}(x;r,s) = \prod_{j=1}^{n}[x+(j-1)(r-s)]^{\underline s}$$
Then  we get  (2.8)  from [10]  i.e.
$$[0,1,2,...,k;b_{ns}(...;r,s)]= \frac{1}{k!}
\sum_{l=s}^{k} {(-1)^{k-l}b_ns(l;r,s)}\Big({k \atop l}\Big)\equiv
S_{r,s}(n,k)$$ becomes the definition of the generalized Stirling
numbers (see A.3.II. - below), which appear to be special case of
general Newton-Stirling numbers of the second kind. (Here
$b_{ns}(...;r,s)(x)= b_{ns}(x;r,s)$.) Naturally the Dobinski-like
formula (2.1) from [10] for exponential polynomials determined by
$[0,1,2,...,k;b_{ns}(.;r,s]= S_{r,s}(n,k)$  is special  case of
(N-S-Dob)  Dobinski-like formula with counting adapted to the
choice $f = b_{ns}$. Along with Bell numbers` sequence or Bessel
numbers`s sequence this special case of Newton-Bell numbers`
sequence$$ B_{r,s}(n)=\sum_{l=s}^{ns} S_{r,s}(n,k)$$  is a moment
sequence [43].\\
\textbf{Example II} The next example of Newton-Stirling numbers
$d_{n,k}$ comes from the paper [44] on  interpolation series
related to the Abel-Goncharov problem. There the divided
difference functional $\Delta_k$ is applied to $e_n$ yielding
$d_{n,k}$ accordingly:
$$\Delta_ke_n=[0,\frac{1}{k},\frac{2}{k},...,\frac{k-1}{k},1;e_n]=d_{n,k}.$$
The general rules for Newton-Stirling arrays  allow us to notice
that
$$d_{n,k}=[0,\frac{1}{k},\frac{2}{k},...,\frac{k-1}{k},1;e_n] =\frac
{k^k}{k!}\sum_{r=0}^{k}(-1)^{k-r}\binom{k}{r} \frac {r^n}{k^n} ;
\quad n\geq k \geq 0,$$
hence for corresponding exponential
polynomials we have

$$\varphi_n(x)=\sum_{k=0}^{n}\frac
{k^k}{k!} \sum_{r=0}^{k}(-1)^{k-r}\binom{k}{r} \frac {r^n}{k^n}
x^k,$$ in accordance with the fact [44] that $k^{n-k} d_{n,k}=
\Big\{{n \atop k}\Big\}.$ Derivation of the Dobinski-like formula
we leave as an exercise.\\
\textbf{On $\psi$-extension.} A $\psi$-extension of the above
Touchard`s symbolic definition of exponential polynomials would
start with the defining formula
$$\varphi_{n+1}= x(\varphi +_{\psi} 1)^n,\quad \varphi^{\underline k}= x^k.\quad\quad\quad\quad\quad(\psi-Exp-Pol)$$
resulting in analogous umbra-symbolic identities and with
corresponding Dobinski-like formula as (35) below, where $b_n =
e_n$. Compare these  with (10) from where we have for this case of
$b_n(x)= e_n(x)= x^n, n\geq 0$ the Newton interpolation formula
$$x^n=\sum_{k=0}^{n}[0,1_{\psi},2_{\psi},...,k_{\psi};e_n]
\psi_{\underline k}(x),\quad n\geq0.$$ For the meaning of the
$\psi$-shift  "$+_{\psi}$"  see [1, 2, 26, 31, 33]. This we shall
develop elsewhere. Meanwhile let us continue the limited review of other extensions.\\
\textbf{Plethystic Stirling numbers` extension} The above umbral
extensions as well as the other extensions to be mentioned in what
follows are to be confronted with inventions of plethystic
exponential polynomials, plethystic Stirling numbers of the second
kind and plethystic extension of Bell numbers from [45] which
constitutes an advanced and profound way to reach the merit of the
finite operator calculus representations - this time realized with
vector space of polynomials in the \textit{infinite }sequence of
variables. In [45] Mendez had derived profits from Nava`s
combinatorics of plethysm then developed by Chen to become an
elegant plethystic representation of umbral calculus  so as to
find out also umbral inverses of plethystic exponential
polynomials and related plethystic Stirling numbers of the first
kind (for references see: [45]). The plethystic exponential
polynomials are  then there expressed via Dobinski-like
(plethystic ) formula
 ( see: (35) in [45]) and  the plethystic Stirling numbers of the second
kind (see: (38) in [45]) are expressed via formula extending the
formula for Stirling numbers of the second kind resulting from the
formula for surjections in its after inclusion-exclusion principle
had been applied form. Whether  $\psi$-extension of plethystic
constructs as above is interesting and possible - we leave as an
inquiry for the future. Occasionally note that though Mendez`s
Stirling numbers of the first and second kind are not Whitney
numbers of an appropriate poset they do bear a striking
resemblance to the latter.\\
\textbf{Whitney numbers, statistic, interpretation.} It is well
known [46] that denoting set of $n$ elements partition lattice by
$\Pi_n$ the arrays  $\left[ {{\begin{array}{*{20}c}
{k}\\{r}\end{array}}}\right]$ and $\Big\{{n \atop k}\Big\}$ are
identified (see also  Theorem 1.3 in [18]) as follows
$$\left[ {{\begin{array}{*{20}c}
{n}\\{n-k}\end{array}}}\right]=w_k(\Pi_n)\quad and  \quad\Big\{{n
\atop {n-k}}\Big\}=W_k(\Pi_n)$$
 where $w_k(\Pi_n)$ and $W_k(\Pi_n)$ denote Whitney numbers of the first and second kind
correspondingly. In order to recognize the possible evolvement of
state of affairs while the combinatorics is concerned let us come
back for a while to $q$-extensions realm. There are several
available ways to define combinatorially $\Big\{{n \atop
k}\Big\}_q$  and $\Big\{{n \atop k}\Big\}^\sim_q $ arrays. Most of
these ways are based on on static on set partitions (see for
example [47], [19], [18], [14], [16], [17], [48], [49]). For
example Gessel in [48] gave to $\Big\{{n \atop k}\Big\}^\sim_q $
combinatorial interpretation as generating functions for an
inversion statistics. In another source paper [49] Milne
demonstrated that $\Big\{{n \atop k}\Big\}_q$  may be viewed in
terms of inversions on partitions and that they  count restricted
growth functions using various statistics (see also [16]). We owe
to Milne also the interpretation of  $\Big\{{n \atop k}\Big\}_q$
as sequences of lines in a corresponding vector space over finite
field. In [19] Sagan delivered the major index statistics`
interpretation of $\Big\{{n \atop k}\Big\}_q$  array of
$q$-Stirling numbers of the second kind. After that the authors of
[18] constructed a family $\wp_n(q)$ of posets as $q$-analogues of
the set partition lattice (different from Dowling $q$-analogue) in
such a manner that (Theorem 5.3 in [18]) $$\left[
{{\begin{array}{*{20}c}
{n}\\{n-k}\end{array}}}\right]=w_k(\wp)\quad and  \quad\Big\{{n
\atop {n-k}}\Big\}_q=W_k(\wp)\quad\forall \wp\in \wp_n(q)$$ become
Whithey numbers $w_k(\wp)$ and  $W_k(\wp)$ of the first and second
kind respectively. Whitney numbers for any graded poset may be
looked at as Stirling like numbers. We shall indicate at the end
of this survey a class of substantially new examples of such
Stirling like arrays - after we inform on prefabs` structures.
Meanwhile let us come back to the main challenge of
$\psi$-extensions where we are faced with an
ispiring obstacle.\\
\textbf{Surprise ?} In \cite{1}, \cite{2} a family of  the so
called $\psi$-Poisson processes was introduced i.e. the
corresponding choice of the function sequence $\psi$ leads  to the
Poisson-like $\psi$-Poisson process. Accordingly one would expect
the extension of Dobinski formula to the $\psi$- case - to be
automatic. Of course it makes no  problem to call the numbers
$$B_n(\psi)=\varphi_n(x=1,\psi)$$
the $\psi$ - Bell numbers - whenever it makes sense - for example
either the sequence of these numbers has combinatorial
interpretation and/or the defining series below are convergent:

\begin{equation}\label{eq35}
B_n(\psi)\equiv\varphi_n(x=1,\psi)= exp^{-1}_\psi\sum_{0\leq
k}\frac{k_\psi^n
 }{k_\psi!}.
\end{equation}
The above might be a far reaching generalization of the standard
case [23]. For example - what about the spectacularly natural and
number theoretic important choice: $\psi=
\langle\frac{1}{F_n!}\rangle_{n\geq 0}$ ($\langle F_n \rangle$ -
Fibonacci sequence )? In this connection (Fibonacci binomial
coefficients  [50] are natural numbers!) a question arises whether
one can prescribe eventual arithmetic properties of some of
$\psi$-Bell numbers beyond $q$-extensions to any kind of composite
modules as in [24], [51] or [52] and [53, 54] - see references
therein. The papers just mentioned perform their investigation
mostly in umbra symbolic Blissard language (see the introduction
in [5] for historical remarks on Blissard`s calculus roots). Note
then (see \textbf{0n $\psi$-extension} remark above) that in the
$\psi$-extensions realm one  may formally introduce the
$\psi$-extended umbra symbol $B_{\psi}$ by analogy to the Bell`s
source of the idea article [41] as follows

$$B(\psi)_{n+1}= (B(\psi) +_{\psi} 1)^n,\quad \varphi^{\underline k}= x^k.\quad\quad\quad\quad\quad(\psi-B-umbra)$$
(see symbolic formula (4.7) in [41] p. 264). The above definition
is equivalent to

$$ B(\psi)_{n+1} =
\sum_{k=0}^{n}\binom{n}{l}_{\psi} B(\psi)_k ,\quad n\geq 0.$$ For
the meaning of the $\psi$-shift "$+_{\psi}$" operator - already
implicit in Ward`s paper [55] -   see [1, 2, 26, 31, 33]. See
occasionally substantial reference to Ward [55]  in Wagner`s
article [15] on generalized Stirling and Lah numbers.\\
\textbf{Question.} Summarizing the discussion above - would we
then - beyond the $q$-umbral extensions` realm - would we have
$\psi$-Bell numbers with Poisson - like processes background - and
\textit{not related to} a kind of Stirling numbers extension - at
least in a way we are acquainted with? Or should we introduce
extended Stirling numbers in another way so as to be \textit{not
related to} Poisson - like processes
beyond the $q$-umbral extensions` realm? \\
On this occasion note also  that all $\psi$ extensions of umbral
calculus do not exhaust all possible representations of $GHW$. For
$GHW$ (Graves-Heisenberg-Weyl) algebra the most general
representation of its` defining commutation relation is already
implicit in [56] which serves [1, 2, 33] as the algebraic operator
formulation of Ward`s calculus of sequences [55]. Namely from the
Rodrigues formula (Theorem 4.3. in [56]) with
$$ \hat{x} q_{n-1}(x)= q_n (x)\quad ,\quad Q q_n = nq_{n-1}$$
it follows that
$$ [Q,\hat{x}]= \textbf{1 },\quad \hat{x}= xQ^{`-1}$$
\vspace{2mm} where  $Q$ - a differential operator [56] is a linear
operator lowering degree of any polynomial by one. $Q$ needs not
to be a delta neither $\psi$-delta operator [2, 1, 33]. We deal
with such a case after the choice of admissible sequence [56, 1,
2, 33] different from $\psi$-sequence
$\langle\frac{1}{n!}\rangle_{n\geq 0}$ or
$\langle\frac{1}{(2n)!}\rangle_{n\geq 0}$ in $(B^\sim)$ and
$(L^\sim)$ for $(B^\sim_n(\psi)$ Bell numbers. Then from [57] we
know that basis consisting of the persistent root polynomials
$\langle\psi_{\underline k}(x)\rangle_{k\geq 0}$ does not
correspond to $\psi$-delta operator. However it determines [56] a
differential operator i.e. the linear one lowering degree of any
polynomial by one.  Another possible "rescue" in seeking for the
convenient, efficient structure with natural objects corresponding
to  Stirling or Bell numbers and Dobinski-like formulas in special
cases  are the  exponential structures and prefabs. For example
reading [30] one notices (section 3) that the $q$-analog of the
Stirling numbers of the second kind  description developed by
Morrison (compare with Section 4 in [15] to see in which way it is
complementary) constitutes {\bf the same} example of Ward`ian -
prefab`ian"  extension  as the Bender - Goldman [58] prefab
example to be considered next right now.  As noticed by Morrison
the relevant prefab exponential formula may equally well be
derived
from the corresponding Stanley`s exponential formula in [59].\\
\textbf{Exponential structures versus prefabs. A subcase of Two
General Classes.} Exponential structures and exponential prefabs
are - in Stanley`s words - basically two ways of looking at the
same phenomenon [59]. Before coming over to inspect [59] from the
"Stirling point of view" let us give at first a family of decisive
examples showing that prefabs are all around us in combinatorics
especially when quite free extensions of Stirling numbers are
concerned. The following example contains
such a family.\\
\textbf{ Bender - Goldman - Wagner  Ward - prefab example.} If
corresponding  "prefabian" $\hat q$-Bell numbers  $ B^{pref}_n(
\gamma )$  are defined as sums over $k$ of  $\hat S_q(n,k) $
Stirling numbers of the  $q$-lattice of unordered direct sums
decompositions of the $n$-dimensional vector space $V_{q,n}$ over
$GF(q) \equiv F_q$ in sect. 2 of [15] then   the formula (2.5) in
[15]  shows up equivalent to the Bender-Goldman exponential
formula (17) from [58] - the source paper on prefabs - and in our
$\psi$-extensions` notation formula (17) with $D_n(q)$ from [58]
now reads:
$$
B^{pref}_{\gamma}(x)= \sum_{n\geq
0}B^{pref}_n(\gamma)\frac{x^n}{n_{\gamma}!}=exp\{exp_{\gamma}(x) -
1\}. \quad\quad\quad(\gamma- e.g.f.)
$$
Here
$$ n_{\gamma}! = (q^n - 1)(q^n - q^1)...(q^n - q^{n-1})=
|GL_n(F_q)|,$$  $D_o(q)=1$ by convention  while $D_n(q) \equiv
B^{pref}_n(\gamma)=$ number of unordered direct sums
decompositions of the vector space $V_{q,n}.$ Compare with
formulae (4-6) in [60] representing the completely new class of
combinatorial prefab structures with noncommutative and
$nonassociative$ composition (synthesis) of its objects. Note "The
natural hint" on $\psi$ extensions remark there  right below these
formulae. Coming back to the Bender - Goldman - Wagner  Ward -
prefab example it is to be noticed that this is a special case of
the First Class formula according to the terminology of three
paths of generalizations being developed in [15]. According to us
Wagner justly refers his First class to Ward [55]. We propose to
call this Wagner`s First Class a "Ward`ian - prefab`ian" Class of
extensions as the characterization formula (1.15) in [15]
\textit{after being summed }over $k$ yields exactly  $\psi$-
extension [60] of prefab exponential formula (12) from  [58] where
$\psi= \langle\frac{1}{f_n }\rangle_{n\geq 0}$  in Wagner`s
notation [15].  Note that our notation [1, 2, 33] is consequently
always "Ward`ian". We also advocate by means of the present paper
the attitude of \textbf{Two General Classes}. The Wagner`s  Class
I is in our terms "\textbf{Ward`ian - prefab`ian}" (see $(\gamma-
e.g.f.)$ above) with $F(n,k)$ Ward-prefab  Stirling numbers and
with $\hat S_q(n,k)$ as example. The second general class in our
terms is "\textbf{Newtonian}" and it incorporates Wagner`s Class
II and Class III  with N-W-C Stirling numbers $U(n,k)$ and with
Newtonian $S_q^\sim(n,k)$ and Gould-Carlitz-Milne $S_q(n,k)$ as
examples mutually expressible each by the other one with help of
re-scaling. One may see that really we are dealing here  with the
\textbf{Newtonian} way notifying that our N-W-C formula (20)
extends (1.12) from the Class III of  [15]  and our N-W-C formula
(14)  extends (1.10) from the Class III of [15]. Note also that
(5.28) from [15] via re-scaling coincides with  $q$ - N-W-C case
of (20) i.e. with $(q - Dobinski)$ formula. Here inevitable
questions arise. For the Newtonian General class we have the
extension (20) of Dobinski formula. So what about the corespondent
formula for the Ward`ian - prefab`ian general class?  ... And what
about \textbf{The Two General Ways} of this paper?  The one "easy"
way is Newtonian. The other way seams to contain   the
$q$-extension as a kind of "singular point" in the domain of
$\psi$-umbral extensions. Is there at and beyond  this "singular
point" of the second way an another  non-Newtonian  second path -
all-embracing what was left beyond the first way ? Before an
attempt to answer some of these questions let us encourage
ourselves by just recalling another distinguished example.

\noindent This another crucial "Ward`ian - prefab`ian" example we
owe to Gessel [48] with his $q$-analog of the exponential formula
as expressed by the Theorem 5.2  from [48].\\
We also recall that the $q$-analog of the Stirling numbers of the
second kind investigated by Morrison in Section 3  of  [30]
constitute {\bf the same} example of  Ward`ian  - prefab`ian
extension  as in the Bender - Goldman - Wagner  Ward - prefab
example. As noticed there by Morrison the  $(\gamma - e.g.f.)$
prefab exponential formula may equally well be derived from the
corresponding Stanley`s exponential formula in [59]. Let us then
now come over to these exponential structures of Stanley with an
expected impact on the current considerations  ( for definitions,
theorems etc. see [59]). In this connection we recall quoting
(notation from [59]) an important class of Stanley`s Stirling -
like numbers $\frac{S_{nk}}{M(n)}$ of the second and those of the
first kind Stanley`s Stirling - like numbers
$\frac{s_{nk}}{M(n)}$. Both kinds are characteristic immanent  for
counting of exponential structures (or equivalently -
corresponding exponential prefabs) and inheriting from there their
combinatorial meaning. This is due to the fact [59] that "with
each exponential structure is associated an "exponential formula"
and more generally a "convolution formula"  which is an analogue
of the well known exponential formula of enumerative
combinatorics" [59]. Consequently with each exponential structure
are associated Stirling-like , Bell-like numbers and Dobinski -
like formulas are expected also.

\noindent \textbf{Exponential structures.} Let $\{Q_n\}_{n\geq0}$
be any exponential structure and let  $\{M(n)\}_{n\geq0}$  be its
denominator sequence  i.e.  $M(n) =$ number of minimal elements of
$Q_n$. Let  $|Q_n|$  be the number of elements of the poset $Q_n$

$$ |Q_n| = \sum_{\pi\in Q_n} 1 .$$
Example: For $Q = \langle \Pi_n \rangle_{n\geq1}$ where $\Pi_n$ is
the partition lattice of  $[n]$ we have   $M(n)= 1$.

\noindent Define "Whitney-Stanley" number $S_{n,k}$ to be the
number of $\pi\in Q_n$ of degree  equal to $k\geq 1$ i.e.

$$ S_{n,k} = \sum_{\pi\in Q_n, |\pi|=k} 1 .$$
Define $S_{n,k}$ - generating characteristic polynomials (vide
exponential polynomials) in standard way

$$
W_n(x)= \sum_{\pi\in Q_n} x^{|\pi|}= \sum_{k=1}^{n} S_{n,k} x^k.
$$
Then  the exponential formula ($W_0(x)=1 = M(0)$) becomes

$$ \sum_{n=0}^{\infty}\frac{W_n(x) y^n}{M(n)n!} = exp\{xq^{-1}(y)\}, $$
where
$$q^{-1}(y) = \sum_{n=1}^{\infty}\frac{y^n}{M(n)n!}\equiv exp_{\psi} - 1 ,$$
with the obvious identification of $\psi$-extension choice here.
Hence the polynomial sequence $\langle p_n(x)= \frac
{W_n(x)}{M(n)}\rangle_{n\geq0}$ constitutes the sequence of
binomial polynomials  i.e. the basic sequence of the corresponding
delta operator $\hat Q = q(D)$. We observe then that
$$ p_n(x) = \sum_{k=0}^{n}\frac{S_{n,k}x^k}{M(n)}\equiv\sum_{k=0}^{n}[0,1,2,...,k;b_n]x^k
$$
are just exponential polynomials` sequence  for the equidistant
nodes case i.e. \textbf{Newton-Stirling} numbers of the second
kind $ S^\sim _{n,k}\equiv\frac{S_{n,k}}{M(n)}$. Both numbers and
the exponential sequence are being bi-univocally  determined by
the exponential structure $Q$. This is a special case of the one
already considered and we have as in this "Lupas case" the
Newton-Stirling-Dobinski formula:
$$ p_n(x) = \frac {1}{\exp(x)}\sum_{k=0}^{\infty}\frac{b_n(k)
x^k}{k!}=
\sum_{k=0}^{n}[0,1,2,...,k;b_n]x^k,\quad\quad\quad\quad(N-S-Dob)$$
where $\langle b_n \rangle_{n\geq0}$ is defined by

$$b_n(x)=\sum_{k=0}^{n}S^\sim _{n,k}x^{\underline k}.$$
\textbf{Note} the \textbf{identification}
$b_n(x)=\frac{w_n(x)}{M(n)},$ where
$$
w_n(x)= -\sum_{\pi\in Q_n} \mu(\hat 0,\pi)\lambda^{|\pi|}.
$$
$\mu$ is M\"{o}bius function and $\hat 0$ is unique minimal
element adjoined to $Q_n$.

\noindent Corresponding Bell-like numbers are then given by
$$ p_n(1) = \frac {1}{\exp(x)}\sum_{k=0}^{\infty}\frac{b_n(k)}{k!}=
\sum_{k=0}^{n}[0,1,2,...,k;b_n],\quad\quad\quad\quad(N-S-Bell).$$

\noindent Besides those above - in Stanley`s paper [59]  there are
implicitly present also inverse-dual "Whitney-Stanley" numbers
$s_{n,k}$ of the first kind i.e.

$$ s_{n,k} = -\sum_{\pi\in Q_n, |\pi|=k}\mu(\hat 0,\pi)  .$$
On this occasion and to the end of considerations on exponential
structures and Stirling like numbers let us make few remarks.
\noindent $q$-extension of exponential formula applied to
enumeration of permutations by inversions is to be find in
Gessel`s paper [48] (see there Theorem 5.2.) where among others he
naturally arrives at the $q$-Stirling numbers of the first kind
giving to them combinatorial interpretation. \noindent Recent
extensions of the exponential formula in the prefab language [58]
are to be find in [60]. Then \textbf{note}: exponential
structures, prefab exponential structures  (extended ones -
included) i.e. schemas where exponential formula holds-imply the
existence of Stirling like and Bell like numbers. As for the
Dobinski-like formulas one needs binomial or extended binomial
coefficients` convolution as it is the case with $\psi$-extensions
of umbral calculus in its operator form.

\noindent \textbf{Information.} On the basis of [60] the present
author introduces new prefab posets` Whitney numbers in [61]. Two
extreme in a sense constructions are proposed there. Namely the
author of [61] introduced two natural partial orders: one $\leq$
in grading-natural subsets of cobweb`s prefabs sets [60] and in
the second proposal one endows the set sums of the so called
"prefabiants" with such another partial order that one arrives at
Bell-like numbers including Fibonacci triad sequences introduced
by the present author in [62].

\vspace{0.2cm}

\subsection*{A.3}
\textbf{Other Generalizations in brief.} We indicate here
\textsc{three} kinds of extensions of Stirling and Bell numbers -
including those which appear in coherent states` applications in
quantum optics on one side  or in the extended rook theory on the
other side. In the supplement for this brief account to follow on
this topics let us note that apart from applications to extended
coherent states` physics of quantum oscillators or strings  [6 -
11, 63, 64] and related Feymann diagrams` description [65] where
we face the spectacular and inevitable emergence of extended
Stirling and Bell numbers (consult  also [66]) there exists a good
deal of  work done on $discretization$ of space - time [67] and/or
Schrodinger equation using umbral methods [68] and GHW algebra
representations in particular (see: [67, 68] for references).\\

\noindent \textbf{A.3.I. An analog of logarithmic algebra.} An
extension of binomiality property from the algebra of formal power
series to the algebra of  formal Laurent series \quad \emph{and
then beyond}  leading to the Loeb`s [69] iterated logarithmic
algebra -  was realized by Roman [70, 71] with the basic
Logarithmic Binomial Formula at the start. The Logarithmic
Binomial Formula and the iterated logarithmic algebra may be given
their $\psi$-analog including $q$-analog of the Logarithmic
Fib-binomial Formula as shown in [72].

\vspace{2mm} The extension of the iterated logarithmic algebra
from [69] is the logarithmic algebra of Loeb-Rota [73]. This
generalization of the formal Laurent series algebra retains main
features and structure of an umbral calculus. Among others it
allows for logarithmic analog of Appell or Sheffer polynomials.
These give rise to  Stirling- type formulas already in [69]. In
[74] Kholodov has invented  an analog of the logarithmic algebra
from [69] in the shape of an umbral calculus on logarithmic
algebras. Specifically (Example 3.1 in [74]) the basic logarithmic
algebra constructed via Jackson derivative $\partial_q$ gives
rise to the analog of $q$-Stirling formula. The
\textit{\it{mnemonic natural question}} arises:  are similar
constructs performable for $\psi$-derivatives $\partial_{\psi}$ ?
\vspace{2mm}

\noindent \textbf{A.3.II. Milne`s Dobinski formula.} In the
classical umbral calculus represented by the finite operator
calculus of Rota the clue and source example of delta operator is
$$\Delta = exp\{D\} - I .$$
Naturally  such delta operator generates Stirling numbers of the
second kind via

$$ k!\Big\{{n \atop k}\Big\} = \Delta^k x^n|_{x=0} .   $$
Accordingly clue and source example of delta operator of the
$\psi$- umbral calculus would be  $$\Delta_{\psi} =
exp\{\partial_{\psi}\} - I .$$ However already Milne`s
$q$-extension [12] - contrary to the above - does not rely on
Jackson derivative $\partial_q$ and it reads

$$ k_q!\Big\{{n \atop k}\Big\}_q = \Delta_q^k x_q^n|_{x=0} ,   $$
where  $\Delta_q^k $ - th k-th difference operator is defined
inductively so that
$$ \Delta_q^k = (exp\{D\}-q^{k-1}I)(exp\{D\}-q^{k-2}I)(exp\{D\}-q^{k-3}I)...(exp\{D\}-q^0I)  $$
The corresponding $q$-Dobinski formula ( 1.26 in [12])  looks "
$\psi$-familiar"  (see: (35) above):

\begin{equation}\label{eq36}
B_{q,n+1}= exp_q^{-1}\sum_{1\leq k}\frac{k_q^n
 }{(k-1)_q!}.
\end{equation}
The obvious challenge is an eventual application of that type
extension to other umbral calculi  - including the analog of the
logarithmic  algebra from [73] with use of $\Big\{{n\atop
k}\Big\}^\sim_{\psi}$ perhaps.

\vspace{1mm}

\noindent \textbf{A.3.III. Normal ordering accomplishment and
generating functions as coherent states.} While staying with
formal power series algebra or even its subalgebra of polynomials
- still valuable extension have been applied as desired from at
least two points of view: \textbf{a)} statistics , \textbf{b)}
normal ordering task for quantum oscillator and strings.

\vspace{1mm}

\noindent \textbf{a)} As for statistics recall that $q$-Stirling
numbers of the second kind may be treated as generating functions
for various statistics counting  ( see: [15] and references
therein). This type of role has been given by Wachs and White in
[16] to $p,q$- Stirling number $S_{p,q}(n,k)$ which is generating
function for the two different joint distribution set partitions
statistics. Wachs and White in [16] also hade proposed
interpretations of their $p,q$-analogue of Stirling numbers in
terms of rook placements and restricted growth functions. From the
defining recurrence (4) in [16] one sees that $S_{1,q}(n,k)=
\Big\{{n\atop k}\Big\}^\sim_q$.

\noindent \textbf{b)} Similar two parameter $r,s$-Stirling numbers
$S_{r,s}(n,k|q)$  arise in the normal ordering accomplishment  for
the expression $[(a^+)^ra^s]^n$ [6,7] where $a^+$  and $a$  stay
for creation and annihilation operators for $q$-deformed quantum
oscillator which is equivalent to say that  $ aa^+ - qa^+a=1. $
(For example $a=\partial_q ,\quad a^+ = \hat x$). From the
recurrence (50)  in [7] one sees that $S_{r,1}(n,k|q) =\Big\{{n
\atop k}\Big\}_q $. This special case of Gould-Carlitz-Milne
$q$-Stirling numbers $\Big\{{n \atop k}\Big\}_q $ appears in [8].
The method to use it in order to recognize coherent states as
combinatorial objects was invented by Katriel in [6]. The authors
of  [9-11, 63, 64] develop the consistent scheme of applications
of the properties of $S_{r,s}(n,k|q=1)\equiv S_{r,s}(n,k)$. These
include [11] closed-form expressions for $S_{r,s}(n,k)$, recursion
relations, generating functions, Dobinski-type formulas. Recall
that generating functions are identified with special expectation
values in boson coherent states. Recall that $S_{r,s}(n,k)$  in
terminology proposed in this paper are the special case of
Newton-Stirling numbers of the second kind and correspondingly -
Dobinski-type formulas are - see the Coherent States` Example in
A.2. above. A new perspective opens while considering normal
ordering task not only for  quantum oscillator but also for
strings which are  "many, many oscillators". The first steps had
been spectacular accomplished by the authors of [63]. These
authors obtained not only analytical expressions but also a
combinatorial interpretation of the corresponding "very much
extended" Stirling and Bell numbers. Their properties are
interpreted in [63] in terms of specific graphs. At the same time
the authors of [63] consider an invention of a $q$-analog of [63]
to be "an outstanding problem".

\vspace{2mm}

\noindent \textbf{A.3.IV. From Howard via  Hsu and Yu and Shiue to
Remmenl and Wachs extensions. Information in brief.} Howard`s
[75], via  Hsu`s and Yu`s [76] and Shiue`s [77] to Remmenl and
Wachs [78] \textbf{sequence of extensions} starts with
\textit{degenerate} weighted Stirling numbers [75] used later on
by the authors of [76, 77] to propose respectable, unified
approach to generalized Stirling numbers. This sequence of
extensions ends with elaborated extended rook theory [78] with its
generalized Stirling numbers and $(p,q)$-analogues of Hsu and
Shiue extensions. Recall that $(p,q)$-analogues  of Stirling
numbers were introduced by Wachs and White  in [16]. In more
detail. Hsu and Shiue had provided a unified scheme for many
extensions of Stirling numbers of both kinds known before [77].
They introduced corresponding unified extensions under the
notation:
$$
\overline S^1_{n,k}(\alpha,\beta,r)\quad and \quad\overline
S^2_{n,k}(\alpha,\beta,r),
$$
such that  $\overline S^1_{n,k}(1,0,0)= \left[
{{\begin{array}{*{20}c} {n}\\{k}\end{array}}}\right]$ and
$\overline S^2_{n,k}(1,0,0)= \Big\{{n \atop k}\Big\}$ and
$\overline S^1_{n,k}(\alpha,\beta,r)=\overline
S^2_{n,k}(\beta,\alpha, -r).$

\noindent Guided by Wachs and  White ideas from [16] Remmel and
Wachs have defined in [78] two natural $p,q$-analogues of  Hsu and
Shiue extensions $$\overline S^i_{n,k}(\alpha,\beta,r),\quad
i=1,2.$$ For that to do Remmel and Wachs have used what we would
call the $\psi_{p,q}$ admissible sequence  $\psi_{p,q} = \langle
[n]_{p,q} \rangle_{n \geq 0}$  where $[ \gamma ]_{p,q} =
\frac{p^{\gamma} - q^{\gamma}}{p-q}$  which for  $\gamma \in
N\cup\{0\}$  becomes the known [16] extension of Gauss extension
i.e.  $[n]_{p,q} = q^{n-1} + pq^{n-2}+ ... +p^{n-2}q + p^{n-1}$.
Factorials and $\psi_{p,q}$-binomial coefficients are then defined
accordingly naturally (see: A.1. Notation.)

\vspace{2mm}

\noindent \textbf{A.4. Extended umbral calculus and some
corresponding extensions of Stirling and Bell numbers. Further
examples.} In this part of our presentation we just list some
examples at hand  where evidently  $\psi$-extension is behind the
scenario for special $\psi$-admissible sequence choices.

\noindent  \textbf{Example A.4.1.} In Katriel and Kibler`s
celebrated paper [8] on normal ordering for deformed boson
operators and operator-valued deformed Stirling numbers one uses
the following $\psi$-admissible sequence
$$
\langle\ \psi_n^{p,q} = \frac{1}{[n]!}\rangle_{n\geq 0},\quad[n]=
\frac{q^n-p^n}{q-p},
$$
from Wachs and White source paper [16].

\noindent  \textbf{Example A.4.2.} In Schork`s paper [79] on
fermionic relatives of Stirling and Lah numbers one uses the
following $\psi$-admissible sequence
$$
\langle\ \psi_n = \frac{1}{[n]^F!}\rangle_{n\geq 0},\quad[n]^F=
n_{q=-1}= \frac{1-(-1)^n}{2}.
$$

\noindent  \textbf{Example A.4.3.} In Parthasarathy`s paper [80]
on fermionic numbers and their roles in some physical problems one
uses the following $\psi$-admissible sequence
$$
\langle\ \psi_n = \frac{1}{[n]^f!}\rangle_{n\geq 0},\quad[n]^f= =
\frac{1-(-1)^nq^n}{1+q}.
$$
Here in [80] the $q$-fermion numbers emerging from the $q$-fermion
oscillator algebra are used to reproduce the $q$-fermionic
Stirling and Bell numbers. New recurrence relations for the
expansion coefficients in the 'anti-normal ordering' of the
q-fermion operators are derived. Corresponding $\psi$-extended
Dobinski formula (see: (15) in [80]) is derived.

\noindent  \textbf{Example A.4.4.} In the paper [11] on extended
Bell and Stirling numbers from hypergeometric exponentiation one
uses the following $\psi$-admissible sequence
$$
\langle\ \psi_n^L = \frac{1}{n!^{L+1}}\rangle_{n\geq 0}.
$$Here in [11] elements of $\psi^L$ -umbral
calculus are at work. Among others the corresponding
$\psi^L$-extended Dobinski formula (see: (15) in  [11]) is
derived.

\noindent  \textbf{Example A.4.5.} In the paper [81] on
representations of the so called "Monomiality Principle"  with
Sheffer-type polynomials and boson normal ordering - just the
standard $\psi$-admissible sequence choice
$$
\langle\ \psi_n = \frac{1}{n!}\rangle_{n\geq 0}.
$$
of the classical non-extended umbral calculus is naturally
abounding in uncountable  formal series indicators of delta
operators examples; see: a) - g) page 3 in [81]. Note there also
GHW - algebra formula (12). As for what the authors  re-discover
(?) to be the so called  "monomiality principle"  one should note
and compare this with the source paper [56] from 1978 by George
Markowsky on "Differential operators and Theory of Binomial
Enumeration". In particular see the GHW - algebra in spirit
Theorem 4.3 in [56] (see also [1, 33] for more on that).

\vspace{2mm}

\noindent \textbf{A.5. Historical and bibliographical relevant
remarks.} To this end we shall here  list few \textit{peculiar}
relevant remarks of historical and bibliographical character.

\noindent \textbf{A.5.1.Remark} The history of GHW algebra has its
roots not later then since Graves` work [32] "On the principles
which regulate the interchange of symbols in certain symbolic
equations" from (\textbf{1853-1857}). See [82] by O.V. Viskov "On
One Result of George Boole" (in Russian) from 1997.

\noindent \textbf{A.5.2.Remark} Generalizations  given by formulas
(3) and (4) from Cakic and Milovanovic paper [83] for another
extensions of Stirling numbers of the second kind as well as their
related properties are an old result published in d`Ocagne M.
article in \textbf{1887} [84]. Many other later generalizations
(see [83]) are consequence of Chak`s work [85]  and  Toscano
papers [86-88]. The relevant papers of importance (see: [83]) are
those [89-93] and [75].

\noindent \textbf{A.5.3.Remark} The relevance  of Schl\"{o}milch`s
work  [94] from \textbf{1852} to N-W-C Stirling numbers is taken
down here with pleasure. Another interesting paper refereeing
directly to the  original Dobinski`s work [4] and Dobinski`s point
of view is the  Fekete`s paper [95] from 1999.

\vspace{2mm}

\textbf{Acknowledgements} The author is much indebted to the
Referee for valuable  indications how to improve the scope and the
shape of the presentation. \noindent Discussions with participants
of Gian-Carlo Rota Polish Seminar\\
$http://ii.uwb.edu.pl/akk/index.html$  - are appreciated also.


\vspace{2mm}


\begin{thebibliography}{99}
\parskip 0pt


\bibitem{1}
A. K. Kwa\'sniewski {\it Main  theorems of extended finite
operator calculus} Integral Transforms and Special Functions, {\bf
14} No 6 (2003), 499-516.

\bibitem{2}
A. K. Kwa\'sniewski {\it On Simple Characterizations of Sheffer
$psi$-polynomials and Related Propositions of the Calculus of
Sequences} ,Bulletin de la Soc. des Sciences et de Lettres de
Lodz,  {\bf 52}  Ser. Rech. Deform.  36 (2002), 45-65. ArXiv:
math.CO/0312397.


\bibitem{3}
Rota G. C. {\it The number of partitions of a set} Amer. Math.
Monthly {\bf 71}(1964), 498-504.

\bibitem{4}
G. Dobinski {\it Summierung der Reihe S .... f\"{u}r m = 1, 2, 3,
4, 5, ....} Grunert Archiv (Arch. Math. Phys) {\bf 61}(1877),
333-336.


\bibitem{5}
Di Nardo E. e Senato D.  {\it Umbral nature of the Poisson random
variables}, Algebraic Combinatorics and Computer science: a
tribute to Gian-Carlo Rota (eds. H. Crapo, D. Senato)
Springer-Verlag, (2001), 245-266.


\bibitem{6}
J. Katriel {\it Bell numbers and coherent states} Physics Letters
A, {\bf 273} (3) (2000), 159-161.

\bibitem{7}
M. Schork {\it On the combinatorics of normal ordering bosonic
operators and deformations of it}   J. Phys. A: Math. Gen. {\bf
36} (2003), 4651-4665.

\bibitem{8}
J.Katriel, M. Kibler {\it Normal ordering for deformed boson
operators and operator-valued deformed Stirling numbers} J. Phys.
A: Math. Gen.{\bf 25} (1992), 2683-26-91.

\bibitem{9}
P. Blasiak, K.A.Penson , A. I. Solomon {\it Dobinski-type
Relations and Log-normal distribution} J. Phys. A Gen.{\bf 36}
L273 (2003).


\bibitem{10}
P. Blasiak, K.A.Penson , A. I. Solomon {\it The Boson Normal
Ordering problem and Generalized Bell Numbers} Annals of
Combinatorics {\bf 7} (2003), 127-139.


\bibitem{11}
P. Blasiak, K.A.Penson , A. I. Solomon {\it Extended Bell and
Stirling numbers from hypergeometric exponentiation} J.Integ. Seq.
{\bf 4} article 01.1.4 (2001).


\bibitem{12}
S.C. Milne {\it A $q$-analog of restricted growth functions,
Dobinski's equality, and Charlier polynomials} ,Trans. Amer. Math.
Soc. {\bf 245} (1978), 89-118.


\bibitem{13}
I. Newton {\it Philosophiae Naturalis Principia Mathematica},
Liber III, Lemma V,  London (1687).

\bibitem{14}
Carl G. Wagner {\it    Partition Statistics and q-Bell Numbers (q
= -1)}, Journal of Integer Sequences {\bf 7}(2004) Article 04.1.1


\bibitem{15}
Carl G. Wagner {\it Generalized Stirling and Lah numbers} Discrete
Mathematics  {\bf 160} (1996), 199-218.


\bibitem{16}
Wachs, D. White  {\it    p,q-Stirling numbers and set partitions
Statistics },  J. Combin. Theory (A) {\bf 56}(1991), 27-46.

\bibitem{17}
R.Ehrenborg {\it Determinants involving $q$-Stirling numbers},
Advances in Applied Mathematics {\bf 31}(2003), 630-642.

\bibitem{18}
Bennett Curtis, Dempsey Kathy J. , Sagan Bruce E. {\it Partition
Lattice $q$-Analogs Related to $q$-Stirling Numbers } Journal
ofAlgebraic Combinatorics {\bf 03}(3),(1994), 261-283.


\bibitem{19}
Bruce E. Sagan {\it A maj static for set partitions}   European J.
Combin. {\bf 12}, (1991), 69-79.

\bibitem{20}
Warren P. Johnson {\it Some applications of the $q$-exponential
formula} Discrete Mathematics {\bf 157} (1996), 207-225.

\bibitem{21}
J. Cigler {\it A new $q$-Analogue of Stirling Numbers} Sitzunber.
Abt. II {\bf 201}(1992), 97-109.

\bibitem{22}
A. K. Kwa\'sniewski {\it $q$-Poisson,$q$-Dobinski, $q$-Rota and
$q$-coherent states-a second  fortieth anniversary memoir},  Proc.
Jangjeon Math. Soc.{\bf 7} (2), (2004), 95-98.

\bibitem{23}
H.S.Wilf {\it Generatingfunctionology} Boston: Academic Press,
1990

\bibitem{24}
J.  Touchard {\it Nombres Exponentiells et Nombres de
Bernoulli}Canad. J. Math. {bf  8} (1956), 305-320.


\bibitem{25}
H. W. Gould  {\it $q$-Stirling numbers of the first and second
kind} Duke Math. J. {\bf 28} (1961), 281-289.

\bibitem{26}
A. K. Kwa\'sniewski {\it  On extended finite operator calculus of
Rota and quantum groups}  Integral Transforms and Special
Functions {\bf 2} No4 (2001),333-340.

\bibitem{27}
Bernd Voigt {\it A common generalization of binomial coefficients,
Stirling numbers and Gaussian coefficients} Publ. I.R.M.A.
Strasbourg, 1984, 229/S-08   Actes 8e  Seminaire Lotharingien, p.
87-89.  http://www.mat.univie.ac.at/~slc/opapers/s08voigt.pdf


\bibitem{28}
J. Konvalina {\it A Unifed Interpretation of the Binomial
Coefficients, the Stirling Numbers and the Gaussian Coefficients},
The Am. Math. Month. {\bf 107} (2000), 901.

\bibitem{29}
J. P. S. Kung {\it The cycle structure of a linear transformation
over a finite field}, Linear  Algebra Appl. {\bf  36}  (1981),
141-155.

\bibitem{30}
Kent E. Morrison  {\it q-exponential families} The Electronic
Journal of Combinatorics {\bf 11} (2004) , No R36

\bibitem{31}
A. K. Kwa\'sniewski {\it Towards  $\psi$-extension of Finite
Operator Calculus of Rota} Rep. Math. Phys. {\bf 48} No3 (2001),
305-342.

\bibitem{32}
Graves C.  {\it On the principles which regulate the interchange
of symbols in certain symbolic equations} Proc. Royal Irish
Academy {\bf 6} (1853-1857), 144-152.


\bibitem{33}
A. K. Kwa\'sniewski, E. Borak {\it Extended finite operator
calculus - an example of algebraization of analysis}  Central
European Journal of Mathematics {\bf 2} (5) (2005), 767-792.


\bibitem{34}
L. Carlitz {\it On arrays of numbers}  Amer. J. Math  {\bf 54}
(1932), 739-752.


\bibitem{35}
W. Lang {\it On generalizations of the Stirling number triangles},
J.Integer Sequences {\bf 3 } (2000), 00.2.4.

\bibitem{36}
A. K. Kwa\'sniewski {\it Cauchy $\hat{q}_\psi$-identity and
$\hat{q}_\psi)$-Fermat matrix via $\hat{q}_\psi$-muting variables
of $\hat{q}_\psi$-Extended Finite Operator Calculus}\\
arXiv:math.CO/0403107 v1 5 March 2004

\bibitem{37}
A. K. Kwa\'sniewski {\it $\psi$-Pascal and $\hat{q}_\psi$-Pascal
matrices - an accessible factory of one source identities and
resulting applications}, Advanced Studies in Contemporary
Mathematics, {\bf 10} No2 (2005), 111-120.


\bibitem{38}
J. Cigler  {\it Operatormethoden f\"{u}r q-Identäten} Monatsh.
Math. {\bf 88},(1979), 87-105.


\bibitem{39}
J. F. Steffensen {\it Some recent researches in the theory of
statistics and actuarial science} Published for the Institute of
Actuaries by the Cambridge Press (1930), 29-34.

\bibitem{40}
Steffensen, J.F. (1950) Interpolation, Chelsea, New York
(reprinted from 1927)

\bibitem{41}
E.T. Bell      {\it Exponential polynomials} Ann. Math. {\bf 35}
No 2 (1934),  258-277.

\bibitem{42}
A. Lupas {\it Dobinski-type formula for binomial polynomials}
Studia Univ.Babes Bolyai Math. {\bf 33} no 2 ,(1988), 40-44.


\bibitem{43}
S.Yassai {\it Examples of moments problems related with some
combinatorial numbers }CNRS U.M.R. 7599 , ``Probabilités et
Modeles Aléatoires'' , Prépublication numéro: PMA-855 2003-10-20 ,
http://www.proba.jussieu.fr/mathdoc/textes/PMA-855.pdf


\bibitem{44}
K. Dzhaparidze {\it On Interpolation Series Related to the
Abel-Goncharov Problem, with Applications to Arithmetic-Geometric
Mean Relationship and Hellinger Integrals} Indag. Mathem., N.S.,
{\bf 12} (1) (2001),   55-72.

\bibitem{45}
M.Mendez {\it Plethystic exponential polynomials and plethystic
Stirling numbers} , Stud. Appl. Math.{\bf96}No 1 (1996), 1-8.

\bibitem{46}
S.A. Joni ,G. C. Rota, B. Sagan {\it From sets to functions: three
elementary examples} Discrete Mathematics {\bf 37} (1981),
193-2002.

\bibitem{47}
E.Steingrímsson {\it Statistics on Ordered Partitions of Sets}

                        math.chalmers univ.preprint (1999)

http://www.math.chalmers.se/~einar/preprints.html


\bibitem{48}
I.M. Gessel {\it A q-analog of the exponential formula}  Discrete
Math. {\bf 40} (1982),  69-80

\bibitem{49}

S.C. Milne {\it Restricted growth functions, rank row matching of
partition lattices and $q$-Stirling numbers} ,Adv. in Math. {\bf
43} (1982), 173-196.


\bibitem{50}
A. K. Kwa\'sniewski {\it Fibonomial cumulative connection
constants} Bulletin of the ICA {\bf 44} (2005), 81-92. ArXiv:
math.CO/0406006 1 Jun 2004

\bibitem{51}
G.T. Williams  {\it Numbers Generated by the Function $e^(e^x-1)$}
The American Mathematical Monthly, {\bf 52}  No. 6. (1945),
323-327.


\bibitem{52}
Lunnon, W. F., Pleasants, P.A.B. Stephens, N.M. {\it Arithmetic
properties of Bell numbers to a composite modulus. I.} Acta Arith.
{bf 35} no. 1 (1979,) 1-16.


\bibitem{53}
H.W. Becker, John Riordan {\it The arithmetic of Bell and Stirling
numbers}Amer. J. Math. {bf  70} (1948), 385-394.

\bibitem{54}
J. Anne Gertsch, Alain M. Robert {\it Some congruences concerning
the Bell numbers} Bull. Belg. Math. Soc. {\bf 3}  (1996), 467-475

\bibitem{55}
M. Ward: {\em A calculus of sequences}, Amer.J.Math. Vol.58,
1936,pp.255-266

\bibitem{56}
G. Markowsky {\it Differential operators and the theory of
binomial enumeration} J.Math.Anal.Appl. {\bf63} (1978), 145-155.

\bibitem{57}
A. Di Bucchianico, D. Loeb{\it  Sequences of Binomial type with
Persistent Roots}J. Math. Anal. Appl. {\bf 199 } (1996), 39-58.



\bibitem{58}
E. Bender, J. Goldman   {\it Enumerative uses of generating
functions} , Indiana Univ. Math.J. {\bf 20} 1971), 753-765.


\bibitem{59}
R. Stanley   {\it Exponential structures} , Studies in Applied
Math. {\bf 59} (1978), 73-82.


\bibitem{60}
A.K.Kwa\`sniewski {\it Cobweb posets as noncommutative prefabs }
ArXiv: math.CO/0503286   (2005)

\bibitem{61}
A.K.Kwa\`sniewski {\it Prefab posets` Whitney numbers } ArXiv:
math.CO/0510027   3 Oct  (2005)

\bibitem{62}
A. K.  Kwa\'sniewski, {\it Fibonacci-triad sequences} Advan. Stud.
Contemp. Math. {\bf 9} (2) (2004),109-118.


\bibitem{63}
M. A. Mendez, P. Blasiak , K. A. Penson {\it Combinatorial
approach to generalized Bell and Stirling numbers and boson normal
ordering problem } arXiv : quant-ph/0505180  May   2005


\bibitem{64}
A.I. Solomon , P. Blasiak , G. Duchamp , A. Horzela , K.A. Penson
{\it  Combinatorial Physics, Normal Order and Model Feynman
Graphs} Proceedings of Symposium 'Symmetries in Science XIII',
Bregenz, Austria, 2003 arXiv: quant-ph/0310174 v1  29 Oct 2003


\bibitem{65}
Bender, C.M, Brody, D.C, and Meister, BK  {\it Quantum field
theory of partitions} Journal of Mathematical Physics, {\bf  40},
(1999),  3239-3245.

\bibitem{66}
Bender, CM, Brody, DC, and Meister, BK (2000) {\it Combinatorics
and field theory} Twistor Newsletter 45, 36-39.

\bibitem{67}
A Dimakis, F Müller-Hoissen and T Striker  {\it Umbral calculus,
discretization, and quantum mechanics on a lattice } J. Phys. A:
Math. Gen. {\bf 29}, (1996) 6861-6876

\bibitem{68}
D. Levi, P. Tempesta and P. Winternitz {\it Umbral Calculus,
Difference Equations and the Discrete Schroedinger Equation}
J.Math.Phys. {\bf 45} (2004) 4077-4105.


\bibitem{69}
 D.E. Loeb {\it The iterated logarithmic algebra} Adv. Math. {\bf 86} (1991), 155-234.

\bibitem{70}
S. M. Roman {\it The Logarithmic Binomial Formula}, Amer.
Math.Monthly {\bf 99}, (1992), 641-648.

\bibitem{71}
 S. M. Roman {\it The Harmonic Logarithms and the  Binomial
Formula}, J.Comb. Theory, Series A  {\bf 63} (1993), 143-163

\bibitem{72}
A. K. Kwa\'sniewski {\it The logarithmic Fib-binomial formula}
Advan. Stud. Contemp. Math.  {\bf 9} No. 1 (2004), 19-26.


\bibitem{73}
D.E. Loeb and G.-C. Rota  {\it Formal power series of logarithmic
type} Adv. Math. {\bf  75} (1988), 1-118.

\bibitem{74}
A.N. Kholodov {\it  The umbral calculus on logarithmic algebras}
Acta Appl. Math.{\bf  19 } (1990), 55-76.

\bibitem{75}
F.T. Howard, Degenerate weighted Stirling numbers Discrete Math.
{\bf 57} no. 1-2, (1985), 45-58.

\bibitem{76}
Hsu L.C. and Yu H. Q.  Northeast. Math. J. {\bf 13} no. 4, (1997),
399-405.

\bibitem{77}
Hsu,L.C.  and Shiue,P.J.S. {\it A unified approach to generalized
Stirling numbers} Adv.in Appl.Math. {\bf 20} no. 3 (1998),
366-384.

\bibitem{78}
Remmel, JB, Wachs, ML, {it Rook theory, generalized Stirling
numbers and (p,q)-analogues} Am. J. Phys. {\bf 72}(2004), 1397.

\bibitem {79}
M. Schork {\it Fermionic relatives of Stirling and Lah numbers }
J. Phys. A: Math. Gen. {\bf 36}(2003), 10391-10398.

\bibitem {80}
Parthasarathy, R {\it  Fermionic Numbers and Their Roles in Some
Physical Problems } (2004) arXiv.org: quant-ph/0403216    to
appear in Phys.Lett.A


\bibitem{81}
P. Blasiak , G. Dattoli, A. Horzela ,  A. Penson  {\it
Representations of Monomiality Principle with Sheffer-type
Polynomials and Boson Normal Ordering } ArXiv 2 Apr 2005

\bibitem{82}
O.V. Viskov "On  One Result of George Boole" (in Russian) Integral
Transforms  and Special Functions-Bulletin vol. 1 No2 (1997) p.
2-7


\bibitem{83}
N.P. Cakic, G.V. Milovanovic {\it On generalized Stirling numbers
and polynomials}, Math. Balkanica (N.S.) {\bf 18}, (2004),
241-248.


\bibitem{84}
D'Ocagne M. {\it  Sur une classe de nombres remarquables} American
J. Math. {\bf 9}, no. 4, (1887), 353-380.

\bibitem {85}
A. M. Chak  {\it A Class of Polynomials and Generalized Stirling
Numbers} Duke Math.Jour. {bf 23} (1956):45-55.


\bibitem{86}
Letterio Toscano,  {\it L. Toscano Numeri di Stirling
generalizzati, operatori differenziali e polinomi ipergeometrici}
Pontifica Academia Scientarum, Commentationes, {\bf 3}
(1939):721-757

\bibitem{87}
Letterio Toscano,  {\it Sulla iterazione dell'operatore xD }
Rendiconti di Matematica e delle sue applicazioni. 5  (VIII)
(1949): 337-350


\bibitem{88}
L. Toscano {\it Generalizzazioni dei polinomi di Laguerre e dei
polinomi attuariali} Riv. Mat. Univ. Parma. (2) {\bf 11} (1970):
191-226.

bibitem{89} L. Carlitz. "Weighted Stirling Numbers of the First
and Second kind I".The Fibonacci Quarterly 18(1980), 147-162.

\bibitem{90}
L. Carlitz. {\it Weighted Stirling Numbers of the First and Second
kind II.}  The Fibonacci Quarterly {\bf 18} (1980), 242-257.

\bibitem{91}
L. Carlitz. "Degenerate Stirling, Bernoulli and Eulerian Numbers."
Utilitas Math. {\bf 15} (1979), 51-88.


\bibitem {92}
P. N. Shrivastava. {\it  On  Generalized Stirling Numbers and
Polynomials} Riv. Mat. Univ. Parma (2) {\bf 11}  (1970): 233-237.

\bibitem {93}
R.C. Singh Chandel. {\it Generalized Stirling Numbers and
Polynomials} Publ. Inst. Math. N.S. {\bf 22}  (36)1977:43-48.

\bibitem{94}
V. P. Sinha , G. K. Dhawan. {\it On Generalized Stirling Numbers
and Polynomials} Riv. Mat. Univ. Parma (2) {\bf 10 } (1969):
95-100.

\bibitem{95}
O. Schl\"{o}milch {\it Recherches sur les coefficient des facultés
analytiques}    J. Reine Angew. Math. {\bf 44} (1852) 344-355.


\bibitem {96} Antal E. Fekete {\it Apropos Bell and Stirling
Numbers } Crux Mathematicorum with Mathematical Mayhem {\bf 25}
No. 5 (1999), 274-281.

\end{thebibliography}
\end{document}